\documentclass[11pt,leqno]{article}

\usepackage[square,numbers]{natbib}
\usepackage{layouts}
\usepackage{enumerate}
\usepackage{authblk}

\usepackage[left=2.9cm,right=2.9cm,top=2.9cm,bottom=2.9cm]{geometry}
\usepackage[utf8]{inputenc}
\usepackage{tikz}
\usepackage{todonotes}
\usepackage{eufrak}
\usepackage{hyperref}
\usepackage{pb-diagram}
\usepackage{amsthm}
\usepackage{amsmath}
\usepackage{txfonts}
\usepackage{algorithm}
\usepackage{algpseudocode}
\usepackage{caption}

\newcommand{\executeiffilenewer}[3]{%
\ifnum\pdfstrcmp{\pdffilemoddate{#1}}%
{\pdffilemoddate{#2}}>0%
{\immediate\write18{#3}}\fi%
}

\usepackage{graphicx}
\newcommand{\private}[1]{}

 \newtheorem{Theorem}{Theorem}
 \newtheorem{Lemma}{Lemma}
 \newtheorem{Proposition}{Proposition}
 \theoremstyle{definition}
 \newtheorem{Definition}{Definition}
 \newtheorem{Remark}{Remark}
 \newtheorem{Example}{Example}
\newcommand{\kommentar}[1]{}

\let \nc \newcommand
\let \rnc \renewcommand

 \rnc \phi \varphi \rnc \epsilon \varepsilon

\nc {\R}{{\bf R}} \nc {\Z}{{\bf Z}}\nc {\Q}{{\bf
Q}} \nc {\N}{{\bf N}} 

\nc {\bd}{\begin{description}} \nc {\ed}{\end{description}} \nc {\bi}{\begin{itemize}} \nc {\ei}{\end{itemize}}
\nc {\be}{\begin{enumerate}} \nc {\ee}{\end{enumerate}} \nc {\bdm}{\begin{displaymath}} \nc
{\edm}{\end{displaymath}} \nc {\bea}{\begin{eqnarray*}} \nc {\eea}{\end{eqnarray*}} \nc {\baa}{\begin{alignat*}}
\nc {\eaa}{\end{alignat*}} \nc {\bsp}{\begin{split}} \nc {\esp}{\end{split}} \nc {\beq}{\begin{equation}} \nc
{\eeq}{\end{equation}} \nc {\btab}{\begin{tabular}} \nc {\etab}{\end{tabular}} \nc {\ba}{\begin{array}} \nc
{\ea}{\end{array}}

\nc{\Gr}{{\bf Gr}}
\nc{\RP}{{\bf RP}}

\newcommand{\cC}{{\mathcal C}}
\newcommand{\cS}{{\mathcal S}}

\newcommand{\cH}{{\mathcal H}}
\newcommand{\cP}{{\mathcal P}}

\nc{\LOT}[1]{L^2(\Omega^{0+1+2}(T;#1))} \nc{\LO}[2]{L^2(\Omega^{0+1+2}(#1;#2))} \nc{\HT}[1]{{H^{0+1+2}(T;#1)}}
\nc{\cHT}[1]{{\cH^{0+1+2}(T;#1)}} \nc{\stor}{D^2\times S^1}

\nc{\Om}{\Omega} \nc{\dist}{{\rm dist}}

\nc{\homeo}{\approx} \nc{\im}{\operatorname{\rm
im}}\nc{\diag}{\operatorname{\rm diag}} \nc{\Crit}{\operatorname{\rm
Crit}}\nc{\grad}{\operatorname{\rm grad}}
\nc{\Hess}{\operatorname{\rm Hess}}\nc{\Ad}{\operatorname{\rm Ad}}
\rnc{\O}{\operatorname{\rm O}}\nc{\SO}{\operatorname{\rm SO}} \nc{\SU}{\operatorname{\rm SU}}
\nc{\Cl}{\operatorname{\rm Cl}} \nc{\PD}{\operatorname{\rm PD}}
\nc{\ad}{\operatorname{\rm ad}} \nc{\vol}{\operatorname{\rm vol}}
\nc{\hol}{\operatorname{\rm
hol}} \nc{\re}{\operatorname{\rm Re}} \nc{\Id}{\operatorname{\rm
Id}} \nc{\Mas}{\operatorname{\rm Mas}} \nc{\Stab}{\operatorname{\rm
Stab}} \nc{\SF}{\operatorname{\rm SF}} \rnc{\ker}{\operatorname{\rm
ker}} \nc{\tr}{\operatorname{\rm tr}} \nc{\sign}{\operatorname{\rm
sign}} \nc{\Spec}{\operatorname{\rm
Spec}}\nc{\coker}{\operatorname{\rm coker}}
\rnc{\hom}{\operatorname{\rm Hom}} \nc{\ch}{\operatorname{\rm ch}}
\nc{\End}{\operatorname{\rm End}} \nc{\Aut}{\operatorname{\rm Aut}}
\nc{\lat}{(\frac{1}{2}\Z)^2} \nc{\hz}{\tfrac{1}{2}\Z}
\nc{\la}{\langle} \nc{\ra}{\rangle} \nc {\Ra}{\Rightarrow} \nc
{\La}{\Leftarrow} \nc {\lla}{\longleftarrow} \nc
{\nach}{\rightarrow} \nc {\equ}{\Leftrightarrow} \nc
{\gdw}{\Longleftrightarrow} \nc {\lra}{\longrightarrow} \nc
{\Lra}{\Longrightarrow} \nc {\lmt}{\longmapsto} \nc
{\quer}{\overline} \nc {\tensor}{\otimes} \nc {\Rt}{\widetilde\R}
\nc {\contract}{\lrcorner} \nc{\sgn}{\operatorname{\rm sgn}}
\newcommand{\kom}[1]{}

\nc{\lcm}{\operatorname{\rm lcm}}
\nc{\CS}{\operatorname{\rm CS}} \nc{\Ch}{\operatorname{\rm Ch}} \nc{\Td}{\operatorname{\rm
Td}}\nc{\Rank}{\operatorname{\rm Rank}} \nc{\GL}{\operatorname{\rm GL}} 
\nc{\proj}{\operatorname{\rm proj}}

\begin{document}

\title{Geometry of Music Perception}

\author[1]{Benjamin Himpel}
\affil[1]{Reutlingen University\\
	Department of Computer Science\\
	Alteburgstr. 150\\
	72762 Reutlingen}
\renewcommand\Affilfont{\itshape\small}



\date{}




\maketitle

\begin{abstract} 
	
Prevalent neuroscientific theories are combined with acoustic observations from various studies 
to create a consistent geometric model for music perception in order to rationalize, explain and predict psycho-acoustic phenomena. The space of all chords is shown to be a Whitney stratified space. Each stratum is a Riemannian manifold which naturally yields a geodesic distance across strata. The resulting metric is compatible with voice-leading satisfying the triangle inequality. The geometric model allows for rigorous studies of psychoacoustic quantities like roughness and harmonicity as height functions. In order to show how to use the geometric framework in psychoacoustic studies, concepts for the perception of chord resolutions are introduced and analyzed.
\end{abstract}

\section{Introduction}

Jacob Collier's fascinating a cappella arrangement of ``In The Bleak Midwinter''\footnote{\href{https://youtu.be/mPZn4x3uOac}{https://youtu.be/mPZn4x3uOac}, published on December 14, 2016, accessed on July 21, 2022} modulates from the key of E to the key of G half-sharp between the third and fourth verses. This is by design, and he explains this choice in his own metaphorical language\footnote{\href{https://youtu.be/5vrhKI7JHQc?t=946}{https://youtu.be/5vrhKI7JHQc?t=946}, published on December 20, 2018, accessed on July 21, 2022}. In response to the question ``Why does music theory sound good to our ears?'' on Wired.com Tech Support (on May 26, 2021), Jacob Collier answers ``Music theory doesn't really sound like anything. It sounds like parchment. Music sounds like stuff though, and the truth is no one knows. It's a bit of a mystery.''\footnote{\href{https://www.wired.com/video/watch/tech-support-jacob-collier-answers-music-theory-questions-from-twitter}{https://www.wired.com/video/watch/tech-support-jacob-collier-answers-music-theory-questions-from-twitter}, accessed on July 21, 2022} This work addresses precisely the question how to geometrically model what music sounds like. We approach this question like a theoretical physicist would: the world consists of physical objects goverened by differential equations.


\subsection{Background}

Music is based on a temporal sequence of pitched sounds. Over time, theorists have analyzed patterns in musical works and described some classes of tones, sounds and sequences thereof as pitches, chords (harmonies) and melodies/chord progressions, respectively. The resulting theory is used in turn by composers to describe their musical inceptions and allow musicians to reproduce them. The theory of harmonies is also used by jazz musicians as a common basis for spontaneous musical creations. 

There is a lot of research related to our differential-geometric approach to music perception. However, music psychology and music theory remain practically distinct as it was already noted by Carol Krumhansl in 1995 \cite{Krumhansl1995}. She empirically develops in \cite{Krumhansl.2001} a tonal hierarchy in specific musical contexts such as scales and tonal music. Frieder Stolzenberg \cite{Stolzenburg2015Harmonyperceptionperiodicity} presents a formal model for harmony perception based on periodicity detection which is compatible with prior empirical results. Harrison \& Pearce \cite{Harrison2020Simultaneousconsonancemusic} reanalyse and formalize consonance perception data from 4 previous major behavioral studies by way of a computer model written in R. Their conclusion is that simultaneous consonance derives in a large part from three phenomena: interference, periodicity/harmonicity, and cultural familiarity. This suggests that chord pleasantness is a multi-dimensional phenomenon, and experiment design in the study of pleasantness in chord perception is highly problematic. They extend their ideas to introduce a new model for the analysis and generation of voice leadings \cite{Harrison2020ComputationalCognitiveModel}. Marjieh et al. \cite{Marjieh2022} provide a detailed analysis of the relationship between consonance and timbre. A speculative account on the evolutional aspect of consonance has been discussed in \cite{Harrison2021} with the conclusion that understanding evolutionary aspects require elaborate cross-cultural and cross-species studies. Chan et al. \cite{Chan2019ScienceHarmonyPsychophysical} combine the ideas of periodicity and roughness in the language of wave interferences in order to define stationary subharmonic tension (essentially the ratio of a generalization of roughness to different frequencies and periodicity) and use it to develop a new theory of transitional harmony, also known as tension and release. Tonal expectations have been analyzed from a sensoric and cognitive perspective in \cite{Collins2014}. Dmitri Tymoczko \cite{Tymoczko.2006, Tymoczko.2011} provides a geometric model of musical chords. 
He also analyzed three different concepts of musical distance and observed that they are in practise related \cite{Tymoczko.2009}. Since our pitch perception is rather forgiving and imprecise, pitch perception corresponds to a probability distribution and therefore a smoothing should be applied to frequencies as in \cite{MilneAndrewJ..2011}, which gives a rigorous way of evaluating similarity of chords or more generally pitch collections using expectation tensors. Differential geometry has also been used in mathematical musicology by way of gauge theory with the aim of explaining tonal attraction \cite{Graben2017,Blutner2020}. Music was viewed as a dynamical system in order to study tonal relationships \cite{Large2010} or musical performances \cite{Burrows1997}. On the level of audio signals, \cite{Gazor2021,Gazor2022} use Hopf bifurcation control to study sound changes in music. In \cite{Pozo2022, Pozo2022a} music theory for classical and jazz music is formalized by providing a mathematical model for tonality, voice leading and chord progressions, which is very different from the geometric and psychoacoustic approach presented in this paper but could help in further developing it. Recent work by Wall et al. \cite{Wall2020} analyzes voice leading and harmony in the context of musical expectancy which is precisely the motivation for our geometric model. Some very interesting vertical ideas on a scientific approach to music can be found in \cite{wilkerson2014harmony}, even though there are---strictly speaking---no new results in that specific article: The brain's exceptional ability for soft computing and pattern recognition on incomplete or over-determined data is relevant for our model. Microtonal intervals have been discussed in the context of harmony by \cite{Bailes2015,Bridges2008}. Several results from cognitive neuroscience studies in the context of music perception also need to be considered for a geometric model \cite{Leino2007,Zhang2018,PagesPortabella2019,Vuust2022,Sauve2021}.
A more conventional and more elaborate account on a scientific approach to music can be found in \cite{feng2012music}. William Sethares wrote a comprehensive analysis on musical sounds based on roughness \cite{sethares2005tuning}.

\subsection{Aims}

\private{This question can be addressed on various levels, and different hypotheses can be tested. It could be the result of an underlying Markov process in the way we acquire a musical preference. This would mean, that we have just one of many possible outcomes of musical development.\todo[inline]{Add reference.}} We hypothesize that there exists a simple underlying mathematical model and mechanism, which is responsible for the harmonic and melodic development in music, in particular Western music. In order to study changes in sound and time, and since sound and time are best modelled as continuous spaces, we need differential geometry in order to study or construct musical trajectories on these spaces. Since the brain has not been understood well enough, there is currently no way of rigorously proving the correctness of a geometric model by deducing it from the way our brain processes music, even though there is a bit of work in this direction \cite{Tramo.2001,Dumas2013MelodiesInSpace}. Instead, the goal of this research project is to validate the model by verifying its music theoretic implications. Our aim is to provide a framework from a differential geometer's point of view in the spirit of \cite{Rickles.2016, Tymoczko.2011}, which is flexible enough to allow for various existing and forthcoming approaches to studying perceptive aspects of the space of notes and chords. In particular, this will remedy all the limitations of geometric models mentioned in \cite[119ff.]{Krumhansl.2001} by making the relations between notes and chords depend on the context and the order. A focal point for this study is the cadence, ``a melodic or harmonic configuration that creates a sense of resolution'' \cite[p. 105--106]{Randel.1999}, which is an important basis for a lot of modern western music and has a long history in human evolution. It reduces tensions in chords, is related to falling fifths and minimizes voice leading distances. For us, it will serve as a guiding principle for the development of a differential-geometric model. While we hope that this model generalizes to other kinds of music because of its generalist approach, our focus will be on Western music. Note that there are many approaches to analyzing and developing music based on machine learning. For the time being, we will stay away from machine learning, even though we may later use techniques from parameter optimization or artifical neural networks to narrow down the model.

Despite numerous studies on music perception, there is a need for a holistic approach by way of a common computational framework in order to study and compare various psychoacoustic quantities like tension, consonance and roughness in a given context. In the spirit of theoretical physics, we make use of mathematical models, abstractions and generalizations in order to create a geometric framework consistent with prevalent neuroscientific theories and results. We show how to rationalize, explain and predict psychoacoustic phenomena as well as disprove psychoacoustic theories using tools from differential geometry. We will not be able to get as far as explaining Jacob Collier's specific modulation to a half-key, but we will show why half-keys appear naturally from a psychoacoustic point of view. We will describe how this simple yet powerful differential-geometric model opens up new research directions.

\subsection{Main Contribution}

The problem with the abundance of competing approaches to dissonance and tension, apart from the great number of different terminologies, is that they are related but not the same, the neuronal processing behind the perception of music has not been understood and the music theory does not yet have a satisfactory explanation based on existing approaches to dissonance and tension. Our geometric model has been constructed in order for these approaches to be studied, compared, and combined. Despite numerous statistical evaluations of models for dissonance and tension, none of these models can be used directly to compose music or develop music further. The main contribution is therefore to present a new approach to music perception by combining the above approaches to music cognition and geometric modelling in a simple differential-geometric model which can be used together with suitable concepts of consonance and tension to deduce the laws of music theory, lends itself to further research and musical developments, as well as provide a flexible framework to relate the perception of music and music theory. This allows for systematically and quantitatively studying the perception of music and music theory with or without just intonation or various equally or not equally tempered systems and describing new approaches to composition and improvisation in the universal language of mathematics and with the tools provided by geometric analysis. It is general enough and modular so that some or all of the concrete sensoric functions presented here can be replaced with alternative ones. A possible outcome is, that we are able to use certain gradient vectors of psychoacoustic quality functions on the space of chords that explain which chord progressions sound good (at which speed and why) and thereby provide an effective tool for composers.


\subsection{Implications}

Therefore, the aim is not to provide yet another bottom-up approach, but to follow a top-down construction of a convenient model, which integrates roughness, consonance, tension with voice leading in order to be useful for analysing music, composing music and ultimately developing music further. In order to study time-dependent aspects of music, we need to be able to consider derivatives of psychoacoustic functions on a space of musical chords with a Riemannian structure. In particular, we want to associate musical expectation to tension on the space of chords.
Even though many of the underlying ideas can be generalized, we restrict ourselves to Western music for reasons of accessibility and convenience with an octave spanning 12 semitones.

\private{\todo[inline]{There should be a reason, why music developed as it did. Just like the fundamental laws of physics are models for observations in the real world, we want to provide a model for observations in the world of music.}}

\section{Foundations of music perception}

In order to be able to construct a geometric framework which is consistent with the prevalent neuroscientific theories, let us first review and briefly discuss the most relevant results. Human evolution has optimized the ability of our sensory nervous system and our brain to process signals efficiently in order to quickly and easily produce the most useful interpretations and implications. Logarithmic perception of signals \cite{Varshney2013LogarithmicPerception} and pattern recognition \cite{Mattson.2014} are at the heart of this optimized mechanism and also provide the basis for music cognition. Signal detection theory provides a mathematical foundation for constructing psychometric functions as models for music perception \cite{Abdi2010,Wichmann2001,Wichmann2001a,Macmillan2004}. Those readers who are not interested in the underlying mechanisms of music perception are welcome to continue with the mathematical part in Section \ref{sec:pitchspace}.

\subsection{Neural coding}

Sensory organs such as eye, ears, skin, nose, and mouth collect various stimuli for transduction, i.e. the conversion into an action potential, which is then transmitted to the central nervous system \cite{Lodish2000} and processed by the neuronal network in the brain as combination of spike trains \cite{Johnson2000}. Both sensation and perception are based on a physiological process of recognizing patterns in the spike trains \cite{Partridge2003}.

\subsection{Logarithmic perception of signals}\label{sec:logperception}

By the Weber-Fechner law \cite{fechner1860elemente} the perceived intensity $p$ of a signal is logarithmic to the stimulus intensity $S$ above a minimal threshold $S_0$:
\[
p=k\log\left(\frac{S}{S_0}\right).\]
Varshney and Sun \cite{Varshney2013LogarithmicPerception} gave a compelling argument, why this is due to an optimization process in biological evolution where the relative error in information-processing is minimized. Quantization in the brain due to limited resources forces a continuous input signal to be perceived logarithmically. 
The Weber-Fechner law applies to the perception of pressure, temperature, light, time, distance and---most importantly for us---to frequency and amplitude of sound waves.

\subsection{Phase locking}\label{sec:phaselocking}

Synchronization and phase locking is a mechanism in the brain for organizing data, recognizing patterns and soft computing. It has also been proposed and confirmed by Langner in the case of pitches \cite{Langner1997Temporalprocessingpitch, Langner2015NeuralCodePitch}. Phase locking for multiple frequencies has been studied in \cite{Sinz2020}. These pattern recognition capabilities can be explained by human evolution \cite{Mattson.2014}. In \cite[p.193--213]{jordan2013advancing} it is argued, how pattern recognition has improved over millions of years in order to allow for better predictions. It is even suggested that the current age of digitalization adds another layer of neurons to recognize new patterns. Pattern recognition is essential for living beings and humans in particular.

We immediately recognize shapes of objects and rhythmic repetitions of signals. Even if we do not see something clearly, because it is too far away, we can predict the shape within a context and thereby recognize the object. Pattern recognition in signals is based on phase-phase synchronizations. This applies for simultaneously emitted signals like pictures and chords, but also for temporally adjacent patterns like moving pictures and chord progressions. Signal predictions and expectations are based on a continuation of patterns. The more patterns diverge from the predicted patterns the more unexpected a signal is. Arguably, our brain prefers signals where patterns can be detected. Again, possible reasons for this can be found in evolution:
\begin{itemize}
	\item Patterns allow us to predict events, and correctly predicting events allows us to evade dangers or kill pray.
	\item More abstractly, changes of patterns cause a rise of information, and we want to minimize the information we need to process, 
\end{itemize}

According to \cite{Fell2011} processes of our working memory are accomplished by neural operations involving phase-phase synchronization. We can think of working memory as an echo of firing neurons in our brain. Temporally adjacent sounds yield synchronized firings, which not only allow us to detect a rhythm but enable us to detect pitches and relate pitches to each other in chord progressions and melodies.

Quantifying phase synchronizations had been addressed by \cite{Lowet2016} which showed that phase-locking values provide better estimation of oscillatory synchronization than spectral coherence. There are other possible explanations for the relevance of simple ratios and periodicity \ref{sec:periodicity} based on neural coding like cross entropy and minimizing sensoric quantities in the context of estimating distances and other measures. At this point, phase-locking seems to be as good an explanation as any for all kinds of sensatoric phenomena and pattern recognition, even though it will eventually be necessary to confirm this or find better explanation for the signal expectation on the neuronal level. As the mechanism for expectation will be similar for different signals, a geometric model will help to reject explanations and find suitable ones based on psychoacoustic observations.

An example of a popular loss function is cross-entropy which is minimized for the training of artificial neural networks. Since we have a metric on each stratum we can study any height function from a differential geometric point of view. For example we can compute the differential or gradient of the dissonance function by way of which we can find the optimal direction in the space of chords to reduce dissonance as fast as possible. 

Cross-entropy might be a good mathematical concept for the purpose of pattern recognition, where we match information received with the information already stored in the brain.

\subsection{Audio signals}

A vibrating object causes surrounding air molecules to vibrate. As long as the kinetic energy is sustained it spreads as a wave by way of a chain reaction. This sound wave travels through the ear canal into the cochlea. Hair cells inside the cochlea convert the wave into an electrical signal, which then travels along the auditory nerve into the brain.

\begin{figure}[h]
	\centering
	\includegraphics[width=0.7\textwidth]{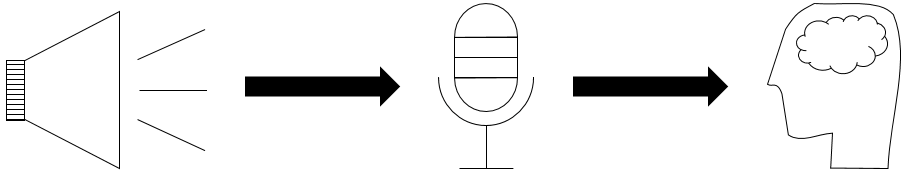}
	\caption{Three major stages of the audio signal's existence \label{fig:sound3ways}}
\end{figure}

The audio signal goes through various stages of existence from the moment of creation to the perception in the brain. Due to a limited resolution of human perception frequency and amplitude is quantized, and the brain logarithmically perceives patterns thereof as certain sound features. These characteristics enable us to quickly recognize and describe instruments, voices and other sounds. We want to distinguish three major stages of an audio signal's existence as shown in Figure \ref{fig:sound3ways}:
\begin{enumerate}
	\item The produced sound, e.g. the vibrating molecules in the air as they are stimulated by a musical instrument or a loudspeaker.
	\item The received sound, e.g. the vibrating microphone diaphragm or the hair cells in the cochlea, at which point the sound wave is converted into an electric signal, before it reaches the brain or different analog or digital recording devices.
	\item The perceived sound, e.g. the interpretation by a person's brain.
\end{enumerate}



\subsection{Spectrum}

The shape of an object is an important factor in the way it can vibrate \cite{Kac1966}. It can be modeled by differential equations involving the geometry of the object. There are several such possibilities known as eigenmodes, each of which moves at a fixed frequency and amplitude as long as the energy is sustained. These eigenmodes are called partials, and the collection of all partials is known as the overtone spectrum of the audio signal. For example, the partials of an ideal vibrating string of length $L$ fixed at both its ends are $n/L$ for $n\in\N$. In this case, the overtone spectrum is called the harmonic spectrum.

Pattern recognition and logarithmic signal perception seem instrumental for the qualitative analysis of sound and music: A musical instrument can play different notes, but our brain detects the same spectral pattern which enables us to identify the sound as coming from the same instrument. This sound quality is also known as timbre and the process of merging several frequencies tonal fusion. Analogous mechanisms apply to voice recognition. Depending on certain deviation patterns in the spectral pattern we can classify and compare different members of the same instrument family (saxophone, clarinet, flute, string, trombone, etc.). It is also exactly this spectral pattern which allows us to recognize the different tones that are played by various sources simultaneously and to determine which instruments are playing which notes, depending on how much training we have.

\subsection{Pitch detection}

Upper partials cannot be easily singled out, only a fundamental frequency can usually be detected by humans. Sounds, where a fundamental frequency can be detected, are called pitched sounds. The process in our brain that detects the pitch is phase locking. The same mechanism is responsible for detecting a pitch in several octaves played together and for detecting a pitch in a tone with a missing fundamental, which seems compatible with autocorrelation \cite{Cariani1996}. Several pitched tones can be played together to produce a chord, where each pitch can be detected. 

Notice that different people might detect different fundamental frequencies depending on the context. This can be seen by considering the ascending Shepard's scale \cite{Shepard.1964} constructed by a series of complex tones which is circular even though the pitch is perceived as only moving upward. 

\subsection{Interference}\label{sec:interference}

Simultaneously emitted Soundwaves interfere with each other. The interference between sine waves with slightly differing frequencies result in beatings which can be computed explicitly. Arbitrary sound waves like those from pitched tones can be approximated by sums of sine waves. The various beatings between slightly different sine wave summands combine to a quality called roughness. Sethares \cite{Sethares1993,sethares2005tuning} uses the Plomp-Levelt curves to provide a formula for measuring roughness and argues that this sound quality is behind tuning and scales. In particular, he suggests that some aspects of music theory can be transferred to compressed and stretched spectra, when played in compressed and stretched scales. This has been confirmed by recent results \cite{Harrison2021CharacterizingSubjectivePleasantness}.

It has been shown by Hinrichsen \cite{Hinrichsen2012} that the tuning of musical instruments such as pianos based on minimizing Shannon entropy of tone spectra is compatible with aural tuning and the Railsback curve. While the tuning of harmonic instruments approximating twelve-tone equal temperament using coinciding partials will work, tuning inharmonic instruments in the context of Western music is more challenging \cite{Cohen1984}. 

Overtone singing is also an interesting aspects of interference. Possibly, overtones are sometimes not what you want to hear, maybe you want to stay away from them, because they are an unwanted artefact.



\subsection{Just-noticeable difference and critical bandwidth}\label{sec:jnd}

The probability for detecting a pitch change between two succeeding tones can be described rigorously using signal detection theory \cite{Abdi2010,Macmillan2004}. It is a collection of psychophysical methods based on statistics for analyzing and determining how signals and noise are perceived.

\begin{figure}[ht]
		\centering
	\input{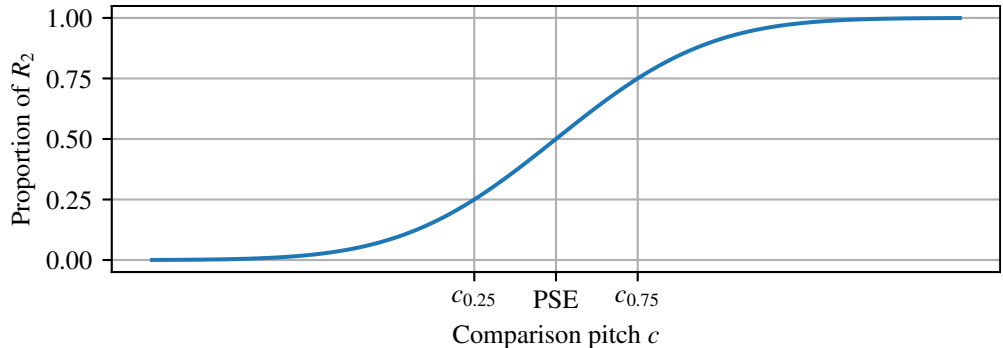}
\caption{Psychometric function with quartiles $c_{0.25}$, PSE and $c_{0.75}$.}\label{fig:jne}
\end{figure}

The just-noticeable difference (JND) also known as difference limen is often described as the minimal difference between two stimuli that can be noticed half of the time. Let us adapt the concise definition and method of computation from psychometric function analysis provided by \cite{Bausenhart2018} to pitch changes. Suppose a subject is presented two succeeding tones as part of a pitch discrimination task. One of the tones is called the reference pitch $p$, the other the comparison pitch $c$. Responses $R_1$ and $R_2$ correspond to the choices $c<p$ and $c>p$, respectively. There is no option $c=p$.  A small set of tone pairs are repeated a number of times (15 to 20), and the subject has to choose one of the two responses. A psychometric function models the proportion of either $R_1$ or $R_2$. For a fixed reference pitch $p$ the psychometric function for $R_2$ should be a monotonically increasing function in the comparison pitch $c$ with values between 0 and 1, because for $c$ much bigger than $p$ the correct response $R_2$ should be obvious. We will assume for simplicity that the shape of the curve fitted to the data follows a cumulative Gaussian as in Figure \ref{fig:jne}, even though other functions like sigmoid, Weibull, logistic or Gumbel are also a possibility \cite{Gilchrist2005}. The point of subjective equality (PSE) is the comparison pitch at which the two responses in this discrimination task are equally likely, i.e. the median. Then the JND is defined to be half its interquartile range, i.e.
\[
{\rm JND} = \frac{c_{0.75}-c_{0.25}}{2},
\]
where $c_{0.25}$ and $c_{0.75}$ represent the comparison pitches, at which a change is detected with probability 0.25 and 0.75, respectively.
\captionsetup[figure]{skip=0pt}
\captionsetup[figure]{skip=10pt}

Notice that when two tones are played in succession the JND is bigger than when the two notes are played simultaneously. This is due to the interference discussed in Section \ref{sec:interference}. Astonishingly, \cite[Section 7.2.2]{HugoFastl2006} states that the JND for two succeeding tones with a pause (difference) is three times higher than the without a pause (modulation). \cite[Figure 7.2]{Rossing2002} shows that the just-noticeable frequency modulation is approximately 3 Hz below 500 Hz and 0.7\% of the frequency above 500 Hz. Clearly, the JND depends on the observer as well as other circumstances (noise) that might interfere with the perception of the signal. 


The critical band is the frequency bandwidth within which the interference between two tones is perceived as beats or roughness, not as two separate tones. The JND is a lot smaller than the critical bandwidth. According to \cite{Smith1999} ``a critical band is 100 Hz wide for center
frequencies below 500 Hz, and 20\% of the center frequency above 500 Hz''. 
A comparison between the critical band and the JND can be seen in \cite[Figure 7.2]{Rossing2002} which in turn is based on \cite[Figure 12]{Zwicker1957}.



In the context of periodicity, Stolzenburg \cite{Stolzenburg2015Harmonyperceptionperiodicity} uses the JND of 1\% and 1.1\% or, equivalently, $\log_2(1.01)\cdot 12 = 0.014355\cdot 12 \approx 17.23$ cent and $\log_2(1.011)\cdot 12 = 0.015783\cdot 12 \approx 18.94$ cent. In \cite{MilneAndrewJ..2011} a standard deviation of 3 cent has been used due to experimentally obtained frequency difference limens of supposedly 3 cent \cite{Moore1984Frequencyintensitydifference}, even though the value of 1\% in \cite{Moore1984Frequencyintensitydifference} corresponds to about 18 cent as we have just seen. Still, the fact that they used the standard deviation of 3 cent for the Gaussian smoothing is an interesting aspect that we will revisit in Section \ref{sec:periodicity}. It will be necessary to design experiments and perform further studies along the lines of \cite{Bubic2009} to collect data for periodicity discrimination in the light of pitch and roughness correlations for tones within chords and between different chords, determine the best model and describe the dependency on noise \cite{Wichmann2001,Wichmann2001a,Hartmann1985}, which is beyond the scope of this work. Due to a lack of such a study, we will assume that cultural familiarity lets us associate slightly mistuned pitches with an ideal pitch and thereby detect and use the implied pitch for the perception of music.


\subsection{Music Perception}\label{sec:music}


Let us define music to be a temporal sequence of pitched sounds created by a formal system. Formal systems obey a set of rules for sound and rhythm, which are ultimately based on physics and mathematics respectively. Different cultures developed and are continuing to develop a variety of systems and scales besides the ones used in Western music \cite{Gill2009}, for example Gamelan music \cite{Becker2019,sethares2005tuning}, Arabic music \cite{Boulos2021,Marcus1993}, Turkish music \cite{Akkoc2002} and classical Indian music \cite{Valla2017}. In Western music there are major subsystems like classical and jazz music. Enculturation is an important factor in the listener's musical expectation and perception \cite{Balkwill.1999,Demorest.2009}, but we want to 
focus on a specific prevalent and in some way universal aspect of music, namely {\em pitch} \cite{Owen.2000,Burton.2015}. 
While the space of received sounds lends itself to a mathematical model, e.g. by using the frequencies and amplitudes computed by Fourier analysis, the spaces of produced and perceived sounds can be compared to it. Given a good microphone connected to some recording device and a good understanding of particle physics the space of produced sounds should be more or less the same as the space of received sounds. Our brain transforms sound waves of music by applying additional filters and perceiving pitch, timbre and loudness. There is also a short term memory effect in the brain, which we hypothesize to be responsible for the sense of resolution in certain chord progressions.

The perception of every person is different and can change via training or degradation. Sound and music are therefore very subjective and can be compared to food, in the sense that the chemical content of food corresponds to the Fourier decomposition of a sound, food can be analysed using chemistry just like we can analyse sound using Fourier analysis or harmony theory, different tastes can be analyzed using signal detection theory and can be described using various characteristics like spiciness, sweetness, sourness, temperature etc. just like sounds can be characterized as warm, loud, sweet, rough etc. via a psychoacoustic analysis. In addition there is an after-taste to food, which might influence the characteristics of food-to-come, just like chord progressions need to be viewed within a musical context.

Chords are also called harmonies and play a key role in Western music. These can sound consonant or dissonant, and the change in this characteristic is an important aspect of musical pieces. Composers build up tension and resolve it subsequently by way of cadences. Notice that it clearly is not only a question of how consonant or dissonant chords sound in a chord progression: the precise way or direction of chord movement is important. It is this kind of aspect in music, that we want to illuminate by geometrically modeling the perception of chords. To this end, we revisit the geometric model of chords \cite{Tymoczko.2006} with a focus on music perception.

\subsection{Mathematics and music}

While sound seems to be well-understood by physics and mathematical structures can be found at every point in music, neither one gives a deep understanding by providing a general principle of how music is perceived by humans. On the other hand, music itself is in reality a mathematical concept based on the brain's perception of sound, put into action in a creative and aesthetically pleasing way: Any kind of scale has been developed mathematically to be compatible with some acoustic observations, rhythm is a time-dependent structure governed by elementary mathematics. Western music theory is a formal system consisting of an assortment of rules that have been deduced from various psychoacoustic preferences.  An account of the major aspects surrounding mathematics and music can be found in \cite{Wright2009}. We want to emphasize the difference between two types of mathematical structures:

The first kind consists of superimposed formal systems in order to give music more structure and to make it more interesting. It starts with simple structures like note lengths and bars to organize rhythm. Other examples include composition procedures like the fugue characterized by imitation and counterpoint as well as various special techniques like Kanon, Krebs, Umkehrung. Then there is the twelve-tone technique invented by Arnold Schoenberg \cite{Schonberg.2000}. For some of these structures we assume a twelve-tone equal temperament, which is itself a mathematical structure superimposed on pitched sounds, not accidentally but deliberately based on a second type of mathematical structure.

This second kind is more subtle, originally due to an evolutionary process and a preference for patterns but ultimately caused by psycho-physical mechanisms like phase locking. It captures the structure inherent in music. It covers temporal structures like rhythmic repetitions. Most Western instruments have approximately a harmonic overtone spectrum. Guided by the simultaneous or sequential perception of intervals and chords humans developed scales, instruments and music theory. Already Pythagoras discovered that simple rational relationships between fundamental frequencies correlate with pleasant sounding intervals. The twelve tones in an octave are also the result of simple rational relationships between frequencies, even though the two physical psychoacoustic qualities harmonicity/periodicity and roughness/interference have been shown to be fundamentally different \cite{sethares2005tuning,Harrison2021CharacterizingSubjectivePleasantness}. Music theory is a formal system which captures more subtle perceptional aspects in Western music. It developed over centuries by the efforts of countless musicians and theorists, mainly however due to observed perceptive qualities of chord progressions. 

Concise models of physical observations can be formulated in the universal language of mathematics, whose powerful tools allow us to deduce complex facts from simple ones. Therefore, the goal is to find a simple way of modeling sounds in the context of music perception, from which we can for example deduce good sounding chord progressions independent of functional harmony, create a music theory in other less common music systems as well as ultimately explain the established Western music theory of harmonies.


\section{Riemannian geometry of chords}\label{sec:pitchspace}


Tymoczko \cite{Tymoczko.2006} viewed the space of chords with $n$ notes as an orbifold \cite{Thurston1997}. In \cite{Hughes2022} the orbifold of chords had been generalized from a topological point of view, while we focus on the  geometry. We argue that it is a Riemannian orbifold \cite{Borzellino1992} and show that the space of chords $\cC$ with an arbitrary numbers of notes is a Whitney stratified space \cite{Pflaum2003AnalyticGeometricStudy} endowed with a metric given by the geodesic distance. The metric provides voice leading distance across different strata. Chord progressions can formally be viewed as sections of the (trivial) $\cC$-bundle over the real line. While our motivation is its use for Western music with its twelve-tone equal temperament, it can readily be adopted to other music. For simplicity, the geometric model represents the chords that can be played using a single instrument which can produce musical tones at any frequency (like a violin) but cannot duplicate notes (like a piano).

Pitches and frequencies can formally be identified with integers via $B3=-1$, $C4=0$, $C\sharp4=1$, etc. Therefore unit distance corresponds to a pitch distance of 100 Cent, which is compatible with the musician's perception of distance between musical tones. The identification between frequency and pitch numbers is given by the function 
\begin{align*}{\rm pitch}: \R &\rightarrow \R,\\
f & \mapsto 12 \cdot \log_2(f/f_0),\end{align*} where $f_0 = 261.626 \rm{Hz}$ corresponds to ${\rm pitch}(f_0)=0=C4$. Chords can then be identified with integer tuples $(p_1,\ldots,p_n) \in \Z^n$. Instead, we will identify chords with tuples $(p_1,\ldots,p_n) \in \R^n$ for the following reasons:
\begin{itemize}
\item There are usually minor pitch adjustments to make chords sound ``better''.
\item The fundamental frequency $f_0$ can assume different values.
\item Quarter tones are entirely legitimate.
\item There are other tuning systems.
\item In particular, not even the piano is tuned using twelve-tone-equal temperament but their stretched tuning follows the Railsback curve \cite{Railsback1938ScaleTemperamentas,Hinrichsen2012}.
\item We assume that instruments play pitches and that the perceived pitch is most relevant for our purpose. We do not include the overtone spectrum with all its amplitudes. When it becomes necessary it can easily be introduced.
\end{itemize}

Since chord notes are played simultaneously, the order of pitches $p_i$ in a chord is irrelevant. For example, the dominant seventh chord (0,4,7,11) needs to be identified with (4,0,7,11). 
\begin{Lemma} Let $S_n$ be the finite symmetric group of all bijective functions $\{1,\ldots,n\} \to \{1,\ldots,n\}$.
\begin{enumerate} 
	\item The permutation \begin{align*}s: \qquad \R^n & \to \R^n\\
		(p_1,\ldots,p_n) &\mapsto (p_{s(1)},\ldots,p_{s(n)}).\end{align*} is a left action on $\R^n$.
	\item The relation \[\tilde c_1 \simeq \tilde c_2 :\Leftrightarrow \exists s\in S_n: s(\tilde c_1)=\tilde c_2 \quad  \text{for } \tilde c_j \in \R^n\] is an equivalence relation on $\R^n$.
\end{enumerate}
\end{Lemma}

\begin{proof} Bijective functions of finite sets form a group.
\begin{enumerate}
	\item We compute that
	\begin{align*}
		s_2(s_1(p_1,\ldots,p_n)) &= s_2(p_{s_1(1)},\ldots,p_{s_1(n)}) = (p_{s_2(s_1(1))},\ldots,p_{s_2(s_1(n))})\\
		& (p_{(s_2\circ s_1)(1)},\ldots,p_{(s_2\circ s_1)(n)}) = (s_2\circ s_1)(p_1,\ldots,p_n).
	\end{align*}
	Therefore $S_n$ acts on $\R^n$ from the left.
	\item Clearly, the relation is reflexive since $\tilde c_1 = \tilde c_1$. If $\tilde c_1 \simeq \tilde c_2$, we have $s(\tilde c_1) = \tilde c_2$ for some $s\in S_n$. Since $S_n$ is a group we have $\tilde c_2 = s^{-1} (\tilde c_1)$. Therefore $\tilde c_2 \simeq \tilde c_1$, and symmetry is satisfied. If $\tilde c_1 \simeq \tilde c_2$ and $\tilde c_2 \simeq \tilde c_3$, then $s_1(\tilde c_1) = \tilde c_2$ and $s_2(\tilde c_2) = \tilde c_3$ for some $s_1,s_2\in S_n$. Therefore $(s_2s_1)(\tilde c_1) = \tilde c_2$ and $\tilde c_1 \tilde \simeq \tilde c_2$ so that the relation is transitive.\qedhere
\end{enumerate}
\end{proof}

Then the quotient by the symmetric group action is given by $\R^n/S_n := \R^n/\simeq$, and its elements are written as $[p_1,\ldots,p_n]$. This space $\R^n/S_n$ is known as the $n$--the symmetric power of $\R$ and is an example of an orbifold\footnote{We can also identify notes with the same name but in different octaves before we consider the quotient by $S_n$. Then we get the toroidal orbifold $(\R/12\Z)^n/S_n$ considered by Tymoczko \cite{Tymoczko.2006,Tymoczko.2011} in order to study efficient voice leading. From a mathematical point of view, this orbifold does not behave differently from $\R^n/S_n$, but this model is not suitable for music perception.} \cite{Thurston1997}, a generalization of a manifold which is locally a quotient of a differentiable manifold by a finite group action. More importantly for us, Theorem \ref{Thm:RiemannianOrbifold} shows that it is a Riemannian orbifold \cite{Borzellino1992,Ratcliffe2007,Lange2020,Bettiol2018} and a Riemannian orbit space \cite{Alekseevsky2003RiemannianGeometryOrbit,Michor2008,Huckemann2010,Thanwerdas2022RiemannianStratifiedGeometries}.

\begin{Definition} A Riemannian orbifold is a metric space which is locally isometric to orbit spaces of isometric actions of finite groups on Riemannian manifolds. A Riemannian orbit space is the quotient of a Riemannian manifold by a proper and isometric Lie group action.
\end{Definition}

\begin{Proposition}\label{Prop:Isometries} Consider the $L^p$ metric on Euclidean space $\R^n$. Then the symmetric group $S_n$ acts on $\R^n$ by isometries.
\end{Proposition}

\begin{proof} Let $(p_1,\ldots,p_n),(q_1,\ldots,q_n)\in \R^n$. Then we get for any $s\in S_n$ by commutativity of the sum
\begin{align*}d((p_1,\ldots,p_n),(q_1,\ldots,q_n)) & = \sum_{k=1}^n (q_k^p-p_k^p)^{1/p} = \sum_{k=1}^n(q_{s(k)}^p-p_{s(k)}^p)^{1/p} \\ & = d((p_{s(1)},\ldots,p_{s(n)}),(q_{s(1)},\ldots,q_{s(n)})).\qedhere
\end{align*}
\end{proof}

This yields the following.
\begin{Theorem}\label{Thm:RiemannianOrbifold}
The quotient space $\cS_n := \R^n/S_n$ is a Riemannian orbifold and a Riemannian orbit space.
\end{Theorem}

In order to study chord progressions it is necessary to consider chords of varying size. We need to construct a metric space of chords with an arbitrary number of tones that is useful for describing music. The metric should provide a sensible voice leading distance, in particular for chord progressions of the form $[0,3] \to [0,3,4]$ or $[0,3] \to [0,3,3]$. 
Multiple same pitches as well as transitions between chords with a different number of tones can be dealt with by considering multiple same pitches in a chord only once, just like a piano plays chords. For example, $[0,0,4,7,11]$ is identified with $[0,4,7,11]$.

\begin{Proposition} Consider the set of chords \[\cS := \left(\bigcup_{k=1}^\infty \cS_k\right).\] The relation
\[\ [p_1,p_2,\ldots,p_k] \sim [p_1,p_2,\ldots,p_{k-1}] :\Leftrightarrow [p_1,p_2,\ldots,p_k] \simeq [p_1,p_1,p_2,\ldots,p_{k-1}] \] for all $k=2,\ldots, n$ is an equivalence relation on $\cS$.
\end{Proposition}

\begin{proof} This is an immediate consequence of $\simeq$ being an equivalence relation.
\end{proof}

This allows us to define the space of all chords.

\begin{Definition} Let $\cC := \cS / \sim$ the space of all chords and $\cC_n := \left(\bigcup_{k=1}^n \cS_k\right)/\sim$ the space of chords with at most $n$ pitches. Let $U_n := \{(p_1,\ldots,p_n) \subset \R^n \mid p_i \neq p_j \text{ for } i\neq j\}$. Let $\pi: \R^n \to \cC_n$ be the quotient map $(p_1,\ldots,p_n) \mapsto [p_1,\ldots,p_n]$.
\end{Definition}

\begin{Remark}
Notice that $\cC_n \setminus \cC_{n-1}$ is the set of chords with exactly $n$ different pitches. 
\end{Remark}

\begin{Example} The space $\cC_2$ is the Euclidean plane as shown in Figure \ref{fig:C2}, where the points are identified with their mirror image when reflected across the diagonal, essentially equivalent to the lower (or the upper) triangle of the plane. $\cC_1$ consists of the singular points with respect to this reflection and is the boundary of $\cC_2$.
\begin{figure}[ht]
		\centering
	\def\svgwidth{6cm}
\begingroup%
  \makeatletter%
  \providecommand\color[2][]{%
    \errmessage{(Inkscape) Color is used for the text in Inkscape, but the package 'color.sty' is not loaded}%
    \renewcommand\color[2][]{}%
  }%
  \providecommand\transparent[1]{%
    \errmessage{(Inkscape) Transparency is used (non-zero) for the text in Inkscape, but the package 'transparent.sty' is not loaded}%
    \renewcommand\transparent[1]{}%
  }%
  \providecommand\rotatebox[2]{#2}%
  \newcommand*\fsize{\dimexpr\f@size pt\relax}%
  \newcommand*\lineheight[1]{\fontsize{\fsize}{#1\fsize}\selectfont}%
  \ifx\svgwidth\undefined%
    \setlength{\unitlength}{181.72118483bp}%
    \ifx\svgscale\undefined%
      \relax%
    \else%
      \setlength{\unitlength}{\unitlength * \real{\svgscale}}%
    \fi%
  \else%
    \setlength{\unitlength}{\svgwidth}%
  \fi%
  \global\let\svgwidth\undefined%
  \global\let\svgscale\undefined%
  \makeatother%
  \begin{picture}(1,0.99999995)%
    \lineheight{1}%
    \setlength\tabcolsep{0pt}%
    \put(0,0){\includegraphics[width=\unitlength,page=1]{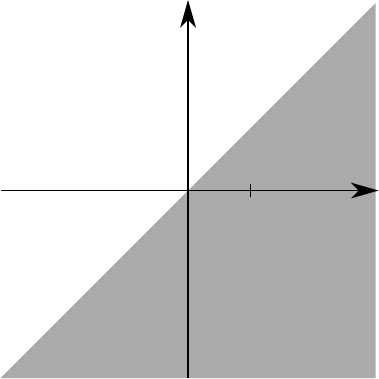}}%
    \put(0.64966033,0.42754734){\color[rgb]{0,0,0}\makebox(0,0)[lt]{\lineheight{1.25}\smash{\begin{tabular}[t]{l}$1$\end{tabular}}}}%
    \put(0.43596516,0.64262737){\color[rgb]{0,0,0}\makebox(0,0)[lt]{\lineheight{1.25}\smash{\begin{tabular}[t]{l}$1$\end{tabular}}}}%
    \put(0,0){\includegraphics[width=\unitlength,page=2]{C2.pdf}}%
    \put(0.59189938,0.85307886){\color[rgb]{0,0,0}\makebox(0,0)[lt]{\lineheight{1.25}\smash{\begin{tabular}[t]{l}$[1,2]$\end{tabular}}}}%
    \put(0.76929056,0.59830783){\color[rgb]{0,0,0}\makebox(0,0)[lt]{\lineheight{1.25}\smash{\begin{tabular}[t]{l}$[2,1]$\end{tabular}}}}%
    \put(0,0){\includegraphics[width=\unitlength,page=3]{C2.pdf}}%
    \put(0.05190114,0.01854){\color[rgb]{0,0,0}\makebox(0,0)[lt]{\lineheight{1.25}\smash{\begin{tabular}[t]{l}$[p]=[p,p]\in \cC_1$\end{tabular}}}}%
  \end{picture}%
\endgroup%

	\caption{The space $\cC_2$ \label{fig:C2}}
\end{figure}
\end{Example}

\begin{Lemma}\label{lem:stab} The Stabilizer $S_p$ of the action of $S_n$ on $\R^n$ is trivial for each $p \in U_n$. 
\end{Lemma}

\begin{proof} The Stabilizer $S_p$ of action of $S_n$ on $\R^n$ is given by $\{s \in S_n \mid s(p) = p\}$. If $s \neq 1$ then $s(i)= j$ for some $i \neq j$. Then $p_i \neq p_{s(j)}$ and therefore $s(p) \neq p$. Therefore $S_p$ is trivial for each $p \in U_n$. 
\end{proof}

\begin{Proposition}\label{prop:stratum} $\cC_n\setminus \cC_{n-1}$ is a Riemannian manifold of dimension $n$. The space of chords $\cC$ is the disjoint union of $\cC_n\setminus \cC_{n=1}$
\[
\cC = \bigsqcup_{n=1}^\infty \left(\cC_n\setminus \cC_{n=1}\right).
\]
\end{Proposition}

\begin{proof} Due to Lemma \ref{lem:stab} we have $[\tilde c_1] \sim [\tilde c_2] \Leftrightarrow \tilde c_1 \simeq \tilde c_2$ for $\tilde c_1, \tilde c_2 \in U_n$. Therefore, $\pi: U_n \to \cC_n \setminus \cC_{n-1}$ is a canonical bijection, and $\cC_n \setminus \cC_{n-1}$ inherits the Riemannian metric from $\R^n$.
\end{proof}

\begin{Remark} The family of chords $\{\cC_n\}_{n\in \N}$ is an example of a filtration 
\[
\cC_1 \subset \cC_2 \subset \ldots \subset \cC_n \subset \ldots
\] By Proposition \ref{prop:stratum} the filtration $\{C_k\}_{k=1,\ldots,n}$ of $\cC$ is an infinite-dimensional stratification, and $\cC_k\setminus \cC_{k-1}$ are the strata of dimensions $k$. 
\end{Remark}

\begin{Remark} The Riemannian metric $g_n$ provides a norm $\|v\|_n$ for every $v \in T_p\cS_n$. Furthermore, the Riemannian metric $g_n$ makes the orbifold $\cS_n$ into a metric space using the geodesic distance defined by 
\[d(p,q) := \inf\left\{\int_{a}^b \|\rho'(t)^p\|_n^{1/p}dt \,\middle\vert\, \begin{split} &\rho: [a,b] \to \cS_n \text{ piecewise smooth},\\ &\rho(a)=p, \rho(b)=q\end{split}\right\} \quad \text{for } p,q \in \cS_n. 
\]
\end{Remark}

\begin{Proposition} The distance on $\cS_n$ can be computed via
\[d_n(p,q)  = \min_{s \in S_n}\tilde d_n(p,s(q))\]
where $\tilde d$ is an $L^p$--metric on $\R^n$.
\end{Proposition}

\begin{proof} In Euclidean space, the geodesic distance is given by the $L^p$--metric. Let $p,q \in \cS_n$. Consider two representatives $(p_1,\ldots,p_n),(q_1,\ldots,q_n)\in U_n$ of $p$ and $q$ with $p_i< p_{i+1}$ und $q_i< q_{i+1}$. Then $tq_i +(1-t)p_i< t q_{i+1} + (1-t)p_{i+1}$ for $t\in[0,1]$ which implies $t(q_1,\ldots,q_n) + (1-t)(p_1,\ldots,p_n) \in U_n$. Therefore $U_n$ is convex. Since $U_n$ is a fundamental domain for $U_n / S_n$, the distance on $U_n / S_n$ is equal to the Euclidean distance in $U_n$. Since the canonical projection $U_n\to \cS_n \setminus \cS_{n-1}$ is an isometric bijection, it follows that for $p,q\in \cS_n \setminus \cS_{n-1}$ we have $d_n(p,q) = \min_{s \in S_n}\tilde d_n(p,s(q))$. Since the closure of $\cS_n \setminus \cS_{n-1}$ is also convex the same formula holds for all $p,q\in \cS_n$.
\end{proof}

Chord progressions in $\cS_n$ with a small distance $d_n$ correspond to efficient voice leading. The metric $d_n$ on $\cS_n$ clearly yields a metric on each stratum $\cC_n\setminus \cC_{n-1}$. Finding a suitable distance on all of $\cC$ is problematic. We can define the following functions $d_n$ and $d$ and $\cC_n$ and $\cC$, respectively:
\[d_n(c_1,c_2)  := \min_{\tilde c_j \in c_j}d_n(\tilde c_1,\tilde c_2) \quad \text{and} \quad
d(c_1,c_2)  := \min_{n\in \N}d_n(c_1,c_2).\]
For example, we compute\begin{align*}d_3([0,1,7],[0,6,7]) & = 5,\\
d_4([0,1,7],[0,6,7]) & = d_4([0,1,7,7],[0,0,6,7]) = 2.\end{align*} 
Even if this is considered to be suitable for determining efficient voice leading, the following shows that this is not a metric on $\cC_n$.

\begin{Proposition} The functions $d_n$ and $d$ on $\cC_n$ and $\cC$ do not satisfy the triangle inequality.
\end{Proposition}
\begin{proof}
Since we have
\begin{equation}d([0],[0,1]) + d([0,1],[0,1,2]) = 1 + 1 < 3 = d([0],[0,1,2]),\label{eq:triangleinequality}
\end{equation}
(see Figure \ref{fig:C3}) this generalization does not satisfy the triangle inequality. The same holds for $d_n$.
\end{proof}

Since the aim is to do differential geometry on $\cC$ the following result is important. See \cite{Pflaum2003AnalyticGeometricStudy} for a detailed treatment of stratified spaces from a geometric analysis point of view.

\begin{Theorem}\label{Thm:Whitney} For each $n \in \N$, the filtration $\{C_k\}_{k\in \N}$ is a Whitney stratification of $\cC$.
\end{Theorem}

\begin{proof} We show that Whitney's condition $B$ is satisfied. Consider the strata $X := C_k \setminus C_{k-1}$ and $Y := C_l \setminus C_{l-1}$ for $k> l$ and embed them in some $\R^N$ via a map $\iota: \cC_k \to \R^N$. Let $x_1,\ldots$ and $y_1,\ldots$ be sequences of points in $X$ and $Y$, respectively, both converging to the same point $y \in Y$, such that the sequence of secant lines $L_i$ between $x_i$ and $y_i$ converges to a line $L \subset \R^N$ in real projective space $\RP^N$ and the sequence of tangent planes $T_i$ to $X$ at the points $x_i$ converges to a $k$--dimensional plane $T$ of $\R^N$ in the Grassmannian $\Gr(k,\R^N)$ as $i$ tends to infinity. The points $x_1,\ldots$ uniquely lift to a sequence $\tilde x_1,\ldots$ in $\R^k$. Let $\tilde y$ be the lift of $Y$ to $\R^k$ so that $\tilde x_1,\ldots$ converges to $\tilde y$. Choose the lift $\tilde Y \subset \R^k$ of $Y$ such that $\tilde y \in \tilde Y$. Then $y_1,\ldots$ uniquely lifts to a sequence $\tilde y_1,\ldots$ of points in $\tilde Y$ that converge to $\tilde y$. Each tangent plane $T_i$ pulls back to the only plane in $\R^k \in \Gr(k,\R^k)$. The secant lines between $(\iota\circ q)^{-1}(x_i)$ and $(\iota\circ q)^{-1}(y_i)$ converge to a line $\tilde L$ in $\RP^k$ which is contained in $\R^k$. This implies that its push-forward $L = d(i\circ q)_{\tilde y} \tilde L$ is contained in $d(i\circ q)_{\tilde y}\R^k = T$.
\end{proof}


Since every stratum of $\cC$ is a metric space and a Riemannian manifold, and the notion of piecewise smooth paths makes sense in $\cC$, we can define the geodesic distance on $\cC$ as follows.

\begin{Definition} We call a continuous path $\rho:[a,b] \to \cC$ piecewise smooth, if there exists a partition $a = x_1 < \ldots x_N = b$ of $[a,b]$ such that $\rho$ restricted to $(x_i,x_{i+1})$ is a smooth path in $\cC_{n_i} \setminus \cC_{n_i-1}$ for some $n_i \in \N$. Let $\rho:[a,b] \to \cC$ be a piecewise smooth path, then we define $\|\rho'(t)\|:=\|\rho'(t)\|_{n}$ if $\rho(t) \in \cC_n \setminus \cC_{n-1}$. The geodesic distance on $\cC$ is \[d(p,q) := \inf\left\{\int_{a}^b \|\rho'(t)^p\|^{1/p}dt \,\middle\vert\, \begin{split} &\rho: [a,b] \to \cC \text{ piecewise smooth},\\ &\rho(a)=p, \rho(b)=q\end{split}\right\} \quad \text{for } p,q \in \cC. 
\]
\end{Definition}

\begin{Theorem}\label{Thm:Distance}
The function $d$ is a metric on $\cC$. It can be computed via
\[
d(p,q) = \inf\left\{\sum_{i=1}^{n-1} d_i(x_i,x_{i+1}) \mid x_i \in \cC_i \setminus \cC_{i-1}\right\}
\]
\end{Theorem}

\begin{proof}
Clearly, $d(p,p)= 0$. Since every stratum is a metric space and we have only a finite number of strata, we get $d(p,q)>0$ for $p\neq q$ and $d(p,q) = d(q,p)$. The concatenation of any piecewise smooth path from $p$ to $q$ and from $q$ to $r$ in $\cC$ is a piecewise smooth path from $p$ to $r$, so that the triangle inequality holds. Therefore, the function $d$ is a metric.

Let $\rho_i:[a,b] \to \cC$ be a sequence of piecewise smooth paths with $\rho_i(a)=p$ and $\rho_i(b)=q$ with a partition $a = x_1 < \ldots x_N = b$ of $[a,b]$ such that $\rho$ restricted to $(x_i,x_{i+1})$ is a smooth path in $\cC_{n_i} \setminus \cC_{n_i-1}$ for some $n_i \in \N$ whose length converges to $d(p,q)$. Since $\cC_{n_i} \setminus \cC_{n_i-1}$ is convex, this implies \[d(p,q) = \sum_{i=1}^{N-1} d_{n_i}(x_i,x_{i+1}).\] Furthermore, we can assume that $n_i > n_{i-1}$ because of this convexity.
\end{proof}

The metric on $\cC$ can be considered as a voice leading distance for music theory. 
\begin{Example}
	Let us compute the distance between $[0]$ and $[0,1,2]$. It can be computed by minimizing the concatenation of geodesic paths within $\cC_3\setminus \cC_2$ and $\cC_2\setminus \cC_1$, and we get \begin{align*}\delta([0],[0,1,2]) & = \min_{p\ge 0} (d_2([0],[0,p]) + d_3([0,p],[0,1,2])) = \min_{p\ge 0} (|p|+|p-1|+|p-2|) \\
		& = 1 + 0 + 1 = 2
	\end{align*}
	In particular, we confirm together with $\delta([0],[0,1])=1$ and $\delta([0,1],[0,1,2])=1$ that the triangle inequality has not been violated as it was in Equation \eqref{eq:triangleinequality}. See Figure \ref{fig:C3}. 
	\begin{figure}[ht]
			\centering
		
		\def\svgwidth{9cm}
\begingroup%
  \makeatletter%
  \providecommand\color[2][]{%
    \errmessage{(Inkscape) Color is used for the text in Inkscape, but the package 'color.sty' is not loaded}%
    \renewcommand\color[2][]{}%
  }%
  \providecommand\transparent[1]{%
    \errmessage{(Inkscape) Transparency is used (non-zero) for the text in Inkscape, but the package 'transparent.sty' is not loaded}%
    \renewcommand\transparent[1]{}%
  }%
  \providecommand\rotatebox[2]{#2}%
  \newcommand*\fsize{\dimexpr\f@size pt\relax}%
  \newcommand*\lineheight[1]{\fontsize{\fsize}{#1\fsize}\selectfont}%
  \ifx\svgwidth\undefined%
    \setlength{\unitlength}{273.08657312bp}%
    \ifx\svgscale\undefined%
      \relax%
    \else%
      \setlength{\unitlength}{\unitlength * \real{\svgscale}}%
    \fi%
  \else%
    \setlength{\unitlength}{\svgwidth}%
  \fi%
  \global\let\svgwidth\undefined%
  \global\let\svgscale\undefined%
  \makeatother%
  \begin{picture}(1,0.60200347)%
    \lineheight{1}%
    \setlength\tabcolsep{0pt}%
    \put(0,0){\includegraphics[width=\unitlength,page=1]{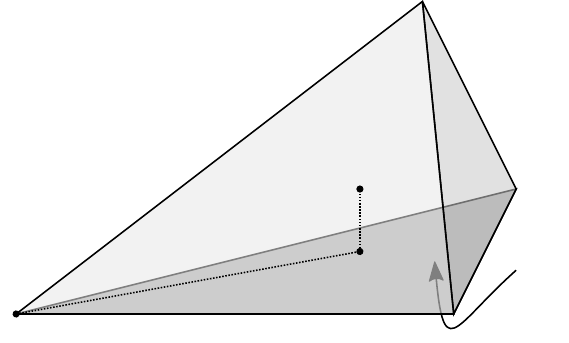}}%
    \put(0.57107078,0.28833254){\color[rgb]{0,0,0}\makebox(0,0)[lt]{\lineheight{1.25}\smash{\begin{tabular}[t]{l}$[0,1,2]$\end{tabular}}}}%
    \put(0.64322872,0.14631312){\color[rgb]{0,0,0}\makebox(0,0)[lt]{\lineheight{1.25}\smash{\begin{tabular}[t]{l}$[0,1]$\end{tabular}}}}%
    \put(-0.00221713,0.00560719){\color[rgb]{0,0,0}\makebox(0,0)[lt]{\lineheight{1.25}\smash{\begin{tabular}[t]{l}$[0]$\end{tabular}}}}%
    \put(0.38416692,0.00879908){\color[rgb]{0,0,0}\makebox(0,0)[lt]{\lineheight{1.25}\smash{\begin{tabular}[t]{l}$\cC_1$\end{tabular}}}}%
    \put(0,0){\includegraphics[width=\unitlength,page=2]{C3.pdf}}%
    \put(0.87851561,0.13787903){\color[rgb]{0,0,0}\makebox(0,0)[lt]{\lineheight{1.25}\smash{\begin{tabular}[t]{l}$\cC_2$\end{tabular}}}}%
  \end{picture}%
\endgroup%

		\caption{The stratum $\cC_3$ \label{fig:C3}}
	\end{figure}
\end{Example}

In summary, 
Theorems \ref{Thm:RiemannianOrbifold}, \ref{Thm:Whitney}, and \ref{Thm:Distance} show that $\cC$ is a well-behaved differential-geometric space:
\begin{enumerate}
	\item $\cC$ is a metric space,
	\item $\cC$ is a Whitney stratified space, and
	\item each stratum of $\cC$ is a Riemannian manifold.
\end{enumerate}
This provides a rich structure for quantitative studies of psychoacoustic models with a voice leading distance on all of $\cC$. The Riemannian metric allows us to study the shape of melodies and chord progressions by differentiating psychoacoustic functions and computing directional derivatives of paths in $\cC$. This model is universal in the sense that it allows note and chord progressions in any musical system. 

\section{Sound perception of chords}\label{sec:heightfunction}

We relate psychometric functions to psychoacoustic height functions on $\cC$. Contour plots of psychoacoustic functions on $\cC$ provide us with insightful visualizations of different models for consonance like roughness and periodicity.\footnote{Timbre and loudness are also important perceptive quantities, which can be addressed later.} The Riemannian structure on $\cC$ allows us to study the shape of melodies and chord progressions as paths in $\cC$ and the perception thereof by differentiating psychoacoustic lifts of the paths in $\cC$ in Section \ref{sec:expectation}.

\subsection{Psychoacoustic functions on the space of chords}


While the space of musical chords can be modelled geometrically, independently of the listener, and a music score can be viewed as a sequence of points or a path in this space, sound perception varies and corresponds to different psychoacoustic functions on this space: dissonance, musical expectation, sense of resolution, root of chord, interference/roughness, all of which depend on both player and listener. Usually, these functions are real-valued on the space of chords (with a given spectrum/timbre) and quantify the individual sensation. This kind of a function turns out to be an example for an important mathematical tool in geometry, analysis and optimization known as a height function on a surface, manifold or more generally a Whitney stratified space. Since the psychoacoustic function varies with the listener and noise it is natural to analyze them using psychometric functions introduced in Section \ref{sec:jnd}. 

Assuming that the JND for pitch discrimination is the same for every reference pitch, let us give a different perspective on the psychometric function from Figure \ref{fig:jne}. Consider the Gaußian distribution $\phi_{\mu,\sigma}$ given by
given by
\[
\phi_{\mu,\sigma}(x) = \frac{1}{\sigma\sqrt{2\pi}}e^{-\frac{1}{2}\left(\frac{x-\mu}{\sigma}\right)^2}.
\]
with mean $\mu = {\rm PSE}$ and standard deviation $\sigma = {\rm JND} / 0.674490$, as well as the Heaviside step function
\[
f(x):= \begin{cases}
	0 & \text{if }  < {\rm 0}\\
	1 & \text{if }  \ge {\rm 0}.
\end{cases}
\]
Then the function in 
Figure \ref{fig:jne} is equal to the convolution $f * \phi_{\mu,\sigma} $ given by\[ (f * \phi_{\mu,\sigma})(p) := \int f(x) \phi_{\mu,\sigma}(p-x)dx .
\]

Consider now a task which is slightly different from the one presented in Section \ref{sec:jnd}: For a given reference pitch $p$ a subject has to say, whether a tone with pitch $c$ has the same pitch as $p$ or not. Let us reformulate the task using random variables. Let $X_{p,c}$ be the random variable which is $1$ (yes) when a comparison pitch $c$ is perceived as the reference pitch $p$ and $0$ (no) otherwise. We can go one step further and consider the continuous random variable $X_p$ which equals $c$ when $p$ is perceived as $c$. Then the probability distribution of $X_{p}$ is given by the normal distribution $\phi_{\mu,\sigma}$ with $\sigma$ and $\mu$ as above. For our purpose let us assume that PSE is equal to $p$.

Again, we can view the probability distribution $\phi_{\mu,\sigma}$ as a convolution of $\phi_{\mu,\sigma}$ with point mass at 0 or, equivalently, as a convolution of $\phi_{0,\sigma}$ with point mass at $\mu=p$.  We observe that $\phi_{\mu,\sigma}(c) = 0.5$ for $c = p\pm{\rm JND}$. Now that we have set up the notation, we can ask which pitch we expect to hear when a tone with pitch $p$ is played. Clearly, it should be $p$, and we can confirm this by computing the expectation value of $X_{p}$:
\[
E(X_{p}) = \int_{-\infty}^\infty P(X_{p}=c) \cdot c \, dc = \int_p^\infty \phi_{\mu,\sigma}(c) \cdot (c + (2p-c)) \, dc = \frac{1}{2}2p = p.
\]
We will use this as a basis for modeling psychoacoustic functions as an expectation value of certain random variables associated to psychometric functions. In \cite{MilneAndrewJ..2011} this viewpoint has been used in order to model perceived distance between pairs of pitch collections, where the perceived dissimilarity was reformulated as a metric between expectation tensors.

As we will see in Section \ref{sec:periodicity}, consonance of dyads and chords is likely to be determined by certain nearby pitches with low periodicity which in turn is due to the phase locking and pattern recognition principle described in Section \ref{sec:phaselocking}. Let us therefore discuss the following multi-variate scenario. Given a fixed set of $N$ pitches $\cP = \{p_1,\ldots,p_N\}$ with $p_i<p_{i+1}$ and ${\rm JND} < |p_{i+1}-p_i| < 2 \cdot {\rm JND}$, a subject has to choose one pitch from $\cP$ which is equal or closest to a given pitch $c$. Let $X_{{\cP},c}$ be the random variable which equals $p_i$ if a perceived pitch $c$ is closest to $p_i$. Clearly, we expect a smoothed version of a step function for the expectation value $E(X_{\cP,c})$ as a function of $p$ where the steps are located at $(p_i+p_{i+1})/2$. One might be tempted to use the convolution of the step function with of $\phi_{\mu,\sigma}$ as above, but by doing so we have neglected the subtle interplay of the random variables and possible dependencies. If we interpret $E(X_{\cP,c})$ as
\[E(X_{\cP,c}) = E\left(\bigvee_{i=1}^N \left(X_{p_i,c}\wedge \bigwedge_{j\neq i}\overline{X_{p_j,c}}\right)\right),\]
we take into account the knowledge that $c$ is not perceived as $p_j$ for $j\neq i$, but we neglect terms of the form $X_{p_i,c} \wedge X_{p_j,c}$ or $\overline{X_{p_1,c}} \wedge \ldots \wedge \overline{X_{p_N,c}}$. If we assume that $X_{p_i,c}$ and $X_{p_j,c}$ are independent random variables for $i\neq j$ we can apply the product formula for independent random variables.

Under the premise that common chord progressions in music theory and their psychoacoustic properties find their justification in certain sound qualities, the chord model $\cC$ together with its sound qualities given by certain height functions on $\cC$ is not only interesting for the music theorist and psychoacoustic analyst, but can become a powerful tool in the hands of composers and computer programs emulating composers because of its conceptual simplicity and quantitative control. Even though $\cC$ could theoretically extend to include the whole overtone spectrum, we hypothesize that different spectra will simply change the psychoacoustic height functions, as long as the spectra consistently have almost the same pattern.

Since there are instruments that do not produce a harmonic series in overtones, it will be interesting to analyse how music and music theory changes for these instruments. A change in the interference scheme due to a different overtone spectrum will promote different note systems. This can be observed in history and other cultures because of the construction of different scales for instruments, which do not produce a harmonic series. Possibly, the relationship of periodicity/harmonicity and consonance needs to be re-evaluated: Is it due to the almost harmonic spectrum of the notes produced by most musical instruments, is it connected to the way human beings interpret periodicity of chords, or are there other more basic concepts at work like logarithmic perception and pattern recognition? However, if it depends on our interpretation of chords, is this due to enculturation or our physical and chemical processing of sound?


In summary, height functions based on mathematical quantitative models for psychoacoustic quantities on the space of chords allow for rigorous studies on music perception. Once the correctness of mathematical models has been confirmed they will yield new music theories. In our work we focus on the psychoacoustic concepts of consonance and tension/release in music. 
From a psychometric point of view it will be necessary to conduct further studies regarding these psychoacoustic quantities. We will see that experiments must be carefully designed as in \cite{Harrison2021CharacterizingSubjectivePleasantness} due to the fortunate (from a Western musical point of view) and at the same time the undesirable (from a scientific point of view) 
correlations between roughness and periodicity.

\subsection{Consonance}

Consonance is a psychoacoustic quality of perceived chords considered to be an important factor in Western music with the usual twelve-tone equal temperament system. Two or more musical tones are considered consonant/dissonant, if they sound pleasant/unpleasant together, and there are a variety of explanations for this phenomenon \cite{Stefano2014UnderstandingMusicalConsonance}. The most important ones go back to roughness (interference) by Helmholtz \cite{Helmholtz1954SensationsTone} and tonal fusion (neural periodicity) by Stumpf \cite{Stumpf2013_1,Stumpf2013_2}. The discussion in \cite{Harrison2020Simultaneousconsonancemusic} carefully analyses various different psychoacoustic interpretations, evaluates data from previous studies, provides a code for several computational models and shows their correlation with consonance ratings. They conclude that consonance depends on interference/roughness, periodicity/harmonicity, and cultural familiarity. While the first two are based on physically justifiable phenomena independent of the individual, cultural familiarity is different for every person by way of musical expertise and cultural conditioning in the following ways:
\begin{enumerate}
	\item Musical training actively and systematically changes your perception. In particular, it allows to better differentiate how consonant chords sound.
	\item The cultural context passively changes your perception by repetition. In particular, it determines how consonant chords sound. E.g. certain jazz chords sound dissonant to people who are unfamiliar with the jazz idiom, while they sound pleasant to jazz musicians.
\end{enumerate} 
Tension, a concept of horizontal harmony between consecutive chords, had also been linked to dissonance \cite{Parncutt2012Consonancedissonancemusic}, but \cite{Lahdelma2020Culturalfamiliaritymusical} suggests that tension is less subjective to cultural familiarity and musical expertise than consonance, pleasantness and harmoniousness of chords. A recent study \cite{Armitage2021} determined that roughness influences automatic responses in a simple cognitive task while harmonicity did not. Furthermore, \cite{Chan2019ScienceHarmonyPsychophysical} argues that tension is independent of harmonicity because it has been shown in \cite{Cook2006PsychophysicsHarmonyPerception} that it is possible for a more consonant chord to resolve into a more dissonant chord. Even though we expect  tension to be related to harmonicity, it is apparently fundamentally different from the vertical quality of consonance and should be reflected in the model accordingly. The difficulty in this discussion surrounding consonance and tension is that in reality they are a conglomeration of different psycho-acoustic phenomena. Furthermore, the terminology might be misleading: Horizontal harmony needs to be viewed in musical context, therefore we will call it the {\em resolve} instead of tension.

Dichotic presentation (different ears for different tones) of chords preserves harmonicity and reduces roughness \cite{Bidelman2009NeuralCorrelatesConsonance}, therefore roughness cannot be responsible for the psychological effect of consonance for chord resolutions, even though roughness and consonance are highly correlated during diotic presentations (same ear for all tones) and will increase the respective effects. The difference of harmonicity and roughness has also been studied in \cite{Tramo.2001}. It is legitimate to say that interference plays a role for the construction of scales, tuning and the quantification of sensory dissonance \cite{sethares2005tuning}, but we hypothesize that there is a fundamental mechanism in the brain that is responsible for the effect of consonance and tension (for a given scale) in the context of chord resolutions and for the way Western music has developed. In particular, such a mechanism should in principal not depend on how badly in tune the notes of a chord sound as long as the chord is approximately correct, and it should not depend on whether the chord tones are presented diotically or dichoticaly. Therefore, we can ignore roughness and beatings for the purpose of studying the mechanism behind chord resolutions. Nevertheless, roughness will strengthen the effect harmonicity has on the listener and will play a role for more subtle variations and fine-tuning of ideal chord progressions.

From a neurophysiological point of view, we hypothesize that roughness, harmonicity and the resolve all find their neural coding origin in the same phase locking principle:
\begin{itemize}
	\item Roughness is based on the interference of sine waves and can be perceived even during dichotic presentation of dyads. It is usually determined using a spectral analysis which will be reviewed briefly in Section \ref{sec:roughness}, but it can also be modeled by the synchronization index model using the degree of phase locking to a particular frequency within the neural pattern \cite{Leman2000} and \cite[Appendix G]{sethares2005tuning}.
	\item Harmonicity can be modeled via periodicity \cite{Stolzenburg2015Harmonyperceptionperiodicity}, which is based on phase locking of perceived pitches and will be discussed in Sections \ref{sec:periodicity} and \ref{sec:chordperiodicity}.
	\item The resolve has not been studied much with respect to the phase locking principle, but we hypothesize that it depends on the interplay between the working memory and harmonicity. Not only has harmonicity been successfully computed via neural periodicity, but working memory has also been linked to phase-phase synchronization \cite{Fell2011}. Some ideas are developed in Section \ref{sec:resolutionsense}.
\end{itemize}

In summary, the three physically justifiable phenomena roughness, harmonicity and resolve are correlated, and their respective psycho-acoustic effects on the listener are amplified by this correlation and by cultural familiarity. These mechanisms are often presented as explanations of consonance, even though they address different issues within the perception of music. The aim of the following sections is therefore to define, distinguish and elaborate upon the individual psycho-acoustic phenomena related to consonance in the context of $\cC$.


\subsection{Roughness}\label{sec:roughness}

Nineteenth-century physicist Herman von Helmholtz \cite{Helmholtz1954SensationsTone} was the first to notice a relation between the harmonic series and the pleasantness of chords, based on which he proposed a theory of consonance and dissonance. In short, he argued, that each tone played by a musical instrument consists of a series of partials determined by the harmonic series: The fewer partials the spectra share, the more dissonant they should be. The interaction between sound waves is called interference, and the interference between two different but similar sine waves create beatings within and roughness outside a critical bandwidth of frequency.

\begin{figure}[ht]
		\centering
	
	\scalebox{0.88}{\input{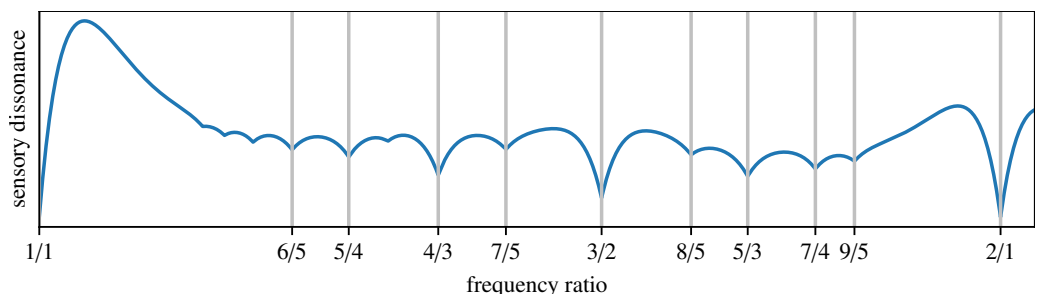}}
	\caption{Sensory dissonance for dyads in terms of their frequency ratio}\label{fig:sensory_dissonance}
\end{figure}

While Western music is usually based on twelve-tone equal temperament, this specific tuning is really a compromise for musical instruments whose pitches are fixed. The pitches of notes for more flexible instruments like the violin or the saxophone are usually adjusted slightly in order to produce chords with minimal or the right amount of roughness. Even pianos are not tuned using twelve-tone-equal temperament but their stretched tuning follows the Railsback curve \cite{Railsback1938ScaleTemperamentas}. Sethares \cite{sethares2005tuning} describes how roughness between complex notes can be computed based on the interference between their partials. He argues that this is one of the main reasons for having a twelve-tone equal temperament system, and that it is a useful tool for tuning and intonating instruments. However, while roughness might be behind tuning, and you want to mostly reduce roughness, it is simply an acoustic artifact that you need to take into account in order to have exactly the correct amount of roughness, just like some coinciding partials whose audibility you want to control. A graph of roughness for dyads can be seen in Figure \ref{fig:sensory_dissonance}\footnote{adapted from https://gist.github.com/endolith/3066664, accessed on October 4, 2021}. The roughness function fits very well into our geometric framework. A contour graph of roughness for triads can be seen in Figure 6.21 of \cite{sethares2005tuning}. One small issue is the fact that the model is not differentiable at its local minima. A possible remedy is the modeling approach by \cite{Leman2000}. Its roughness graph of a harmonic tone complex can be seen in \cite[Figure 4]{Leman2000}. It remains to be seen how deep we have to dive into other aspects like cochlear hydrodynamics \cite{Vencovsky2016} in order to improve the roughness model for further studies on music perception.

\captionsetup[figure]{skip=0pt}
\captionsetup[figure]{skip=10pt}

It will be interesting to study roughness together with the geometric model $\cC$ to compute and visualize entropy and determine ideal tunings. For details, formulas and graphs of roughness we refer to \cite{sethares2005tuning}. More importantly for us we need to study roughness in combination with harmonicity, because both types of consonance are relevant for music in their own way, but their specific psycho-acoustic effects independently of each other are not clear yet.

\subsection{Harmonicity of dyads}\label{sec:periodicity}

In order to create a suitable harmonicity height function $p:\cC\to \R$ on musical chords, we will focus on the harmonicity model determined by relative (logarithmic) periodicity as presented by Stolzenburg \cite{Stolzenburg2015Harmonyperceptionperiodicity}. His explanations based on the neuronal model by Langner \cite{Langner1997Temporalprocessingpitch, Langner2015NeuralCodePitch} using phase locking are convincing, even if probabilistic implications of the psychometric functions apart from the JND have not been included and other aspects like roughness and cultural familiarity clearly alter the perception of chords. In short, if the ratio of two pitch frequencies $f_1$ and $f_2$ with $f_2 \ge f_1$ is given by $f_2/f_1 = p/q$ with $\gcd(p,q)=1$, then the periodicity for this dyad is $q$. In other words, the period of the sound wave for this dyad equals $q$ periods of the first (lower) note. 

Since the above periodicity $q$ will change a lot for small changes of the ratio $f_2/f_1=p/q$, our brain will pick the smallest $q$ within a JND for harmonicity through phase locking as discussed in Section \ref{sec:jnd}. It is chosen to be 1\% and 1.1\% in \cite{Stolzenburg2015Harmonyperceptionperiodicity} based on related results by \cite{Zwicker1957,HugoFastl2006,Hall1984,Roederer2008,Vos1986,Kopiez2003,Moore1984Frequencyintensitydifference,Moore1985,Moore1986,Hartmann2004}. As we have seen in Section \ref{sec:jnd} this corresponds to approximately $18$ cent. 

\begin{figure}[ht]
		\centering
	
	\input{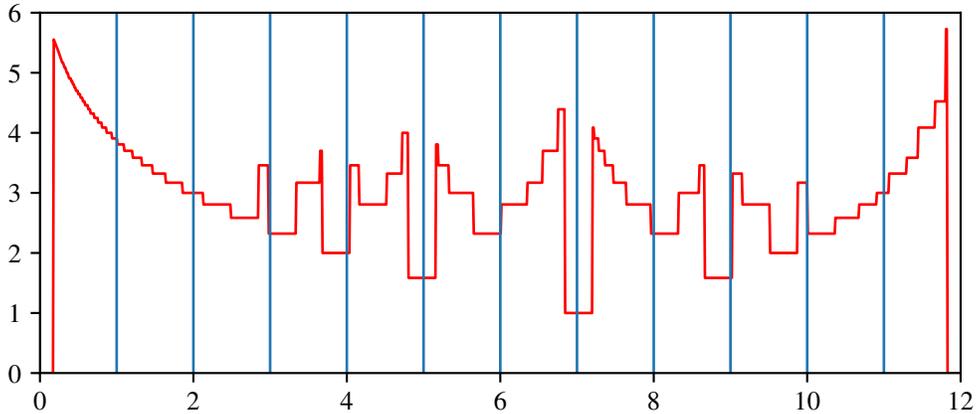}
	\caption{Logarithmic periodicities of dyads spanning at most one octave}\label{fig:periodicity-dyads}
\end{figure}

A naive model for periodicity is therefore given by a step function with a JND of 18 cent for periodicity as shown in Figure \ref{fig:periodicity-dyads}, but we need to keep in mind that depending on the listener, the loudness and distracting noise the JND might vary. Furthermore, in order to compute JND for harmonicity of simultaneously played tones we need to design a new experiment, where we can analyze the effects of roughness and harmonicity separately.

Notice however, that by incorporating probabilistic aspects via Gaussian smoothing after first constructing a step function resembling the periodicity based on \cite{Stolzenburg2015Harmonyperceptionperiodicity} we commit a conceptual error, which we are not able correct in this work but which is hopefully small enough to still provide useful results. The perceived periodicity of a chord is determined by the period that is the best fit for the given spike train induced by the audio signal. The brain either chooses the smallest period it can detect or it detects a mixture of periodicities as an average. It is also possible that different periods are detected at different times within a small time interval due to small variations in the spike sequence or in the pitch. Spike trains with low periodicities are more likely to be detected than spike trains with high periodicities. In order to create a better model we need take into account these probabilistic issues already within the phase locking stage and make use of probabilistic tools like cross entropy and coherence in the time domain along the lines of \cite{Lowet2016}. 

Let us consider a dyad in 12-tone equal temperament as discussed in Section \ref{sec:pitchspace}. If the lower note is fixed, a dyad spanning at most one octave is determined by the number of separating semitones $i\in \{0,\ldots,12\}$. Its frequencies $f_i$ within a JND of 1.1\%, its relative periodicies $L_i$ and its logarithm are given in Table \ref{tab:periodicities}.

\begin{table}[ht]
	\begin{tabular}{c|c|c|c|c|c|c|c|c|c|c|c|c|c}
		$i$ & 0 & 1 & 2 & 3 & 4 & 5 & 6 & 7 & 8 & 9 & 10 & 11 & 12\\
		\hline
		$f_i$ & 1 & $\tfrac{16}{15}$ & $\tfrac{9}{8}$ & $\tfrac{6}{5}$ & $\tfrac{5}{4}$ & $\tfrac{4}{3}$ & $\tfrac{7}{5}$ & $\tfrac{3}{2}$ & $\tfrac{8}{5}$ & $\tfrac{5}{3}$ & $\tfrac{9}{5}$ & $\tfrac{15}{8}$ & $\tfrac{2}{1}$\\
		\hline
		$L_i$ & 1 & 15 & 8 & 5 & 4 & 3 & 5 & 2 & 5 & 3 & 5 & 8 & 1\\
		\hline
		$\log_2(L_i)$ & 0 & 3.91 & 3 & 2.32 & 2 & 1.58 & 2.32 & 1 & 2.32 & 1.58 & 2.32 & 3 & 0
	\end{tabular}
	\caption{Frequencies and periodicites relative to pitch 0}
	\label{tab:periodicities}
\end{table}

\begin{figure}[ht]
		\centering
	
	\input{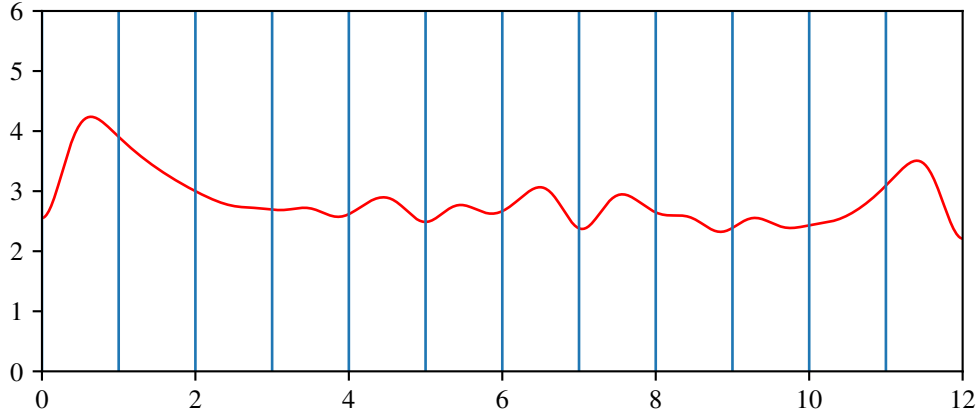}
	\caption{Logarithmic periodicities smoothed by a Gaussian with standard deviation of $\sigma = {\rm JND} / 0.674490 \approx 26.69$ cent}\label{fig:periodicity-dyads-smoothed-largestd}
\end{figure}

Observe that the concept of voice leading is also related to the periodicity of an octave being 1. It allows the player to change the voicing of a chord without changing the psychological effect of its sound by much. Certainly, periodicities can be computed for all intervals as a function $f_p$ for all dyads $[0,p]$, where $p \in [0,12]$. Its graph is shown in Figure \ref{fig:periodicity-dyads}, where the JND is 18 cent. Notice that the step functions has jumps very close to some of the integers. \captionsetup[figure]{skip=0pt}
\captionsetup[figure]{skip=10pt}
As we have discussed above, in order to get a smooth height function on the space of chords in the spirit of psychometric functions we can consider the convolution with a Gaussian. A standard deviation of $\sigma = {\rm JND} / 0.674490$ which we discussed in Section \ref{sec:heightfunction} to be the correct value in the context of psychometric functions seems much too big. When applied to the step function the resulting graph can be seen in Figure \ref{fig:periodicity-dyads-smoothed-largestd}.
\captionsetup[figure]{skip=0pt}
\captionsetup[figure]{skip=10pt}

\begin{figure}[ht]
		\centering
	
	\input{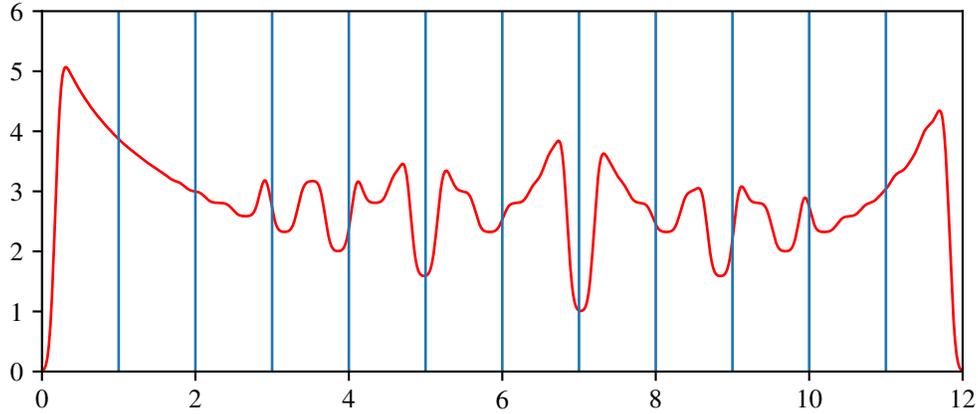}
	\caption{Logarithmic periodicities smoothed by a Gaussian with standard deviation of $\sigma = {\rm JND} / 3 = 6$ cent}\label{fig:periodicity-dyads-smoothed}
\end{figure}

In order to keep the appropriate maxima and minima of the step function a standard deviation of $\sigma = {\rm JND}/3 = 6$ cent seems better. The result is shown in Figure \ref{fig:periodicity-dyads-smoothed}. There are a few reasons why this smaller $\sigma$ is more appropriate. First of all, \cite{MilneAndrewJ..2011} suggests a minimum standard deviation of $3$ cent. Even though this alone is not a good enough reason, especially because \cite{MilneAndrewJ..2011} refers to \cite{Moore1984Frequencyintensitydifference}, where the difference limen has been computed to be approximately 1\%, it suggests that a careful psycho-metric analysis of harmonicity, roughness and pitch needs to be done that sheds some light on their interdependence. We hypothesize that a side effect of roughness is the increase of phase locking precision for the detection of harmonicity. In combination with pitch detection the conditional probability for detecting the correct harmonicity will also increase, because the product of two Gaussians with standard deviations $\sigma_1$ and $\sigma_2$ is again a Gaussian with (smaller) standard deviation \[\sigma = \sqrt{\frac{\sigma_1^2\cdot \sigma_2^2}{\sigma_1^2+\sigma_2^2}}.\] We realize that these reasons need to be elaborated on, treated more rigorously and their effects quantified, but this needs to be done elsewhere. Instead we will lift the periodicity function with its visually and subjectively satisfactory parameters to higher dimensions.

\captionsetup[figure]{skip=0pt}
\captionsetup[figure]{skip=10pt}


\subsection{Harmonicity of arbitrary chords}\label{sec:chordperiodicity}

The definition of periodicity generalizes to chords with more than two notes by letting the periodicity be the smallest positive integer $q$ satisfying $q/f_1=p_2/ f_2= p_n/ f_n$ for some $p_2,\ldots,p_n\in \N$, where $f_1$ is the frequency of the lowest note. Equivalently, periodicity is the smallest positive integer satisfying $f_2/f_1 = p_2/q, f_3/f_1 = p_3/q, \ldots, f_n/f_1 = p_n/q,$ in other words, $q$ is the least common multiple of the denominators in the irreducible fractions representing the frequencies relative to $f_1$. Let $\cC^{[0,12]}_n \subset \cC_n$ be the subspace of all chords where each tone is contained in the octave $[0,12]$ and the base note is equal to 0.\footnote{This can easily be generalized to chords spanning more than an octave.} Define the chords $\cC^p_n \subset \cC^{[0,12]}_n$ with periodicity $p$ via \[\cC^p_n := \left\{\left[0,12\cdot\log_2\left(\tfrac{p_2}{q_2}\right),\ldots,12\cdot\log_2\left(\tfrac{p_n}{q_n}\right)\right] \,\middle\vert\,\\  \begin{split}&\forall_{i} \left(1 \le \tfrac{p_i}{q_i}\le 2 \wedge \gcd(p_i,q_i) = 1\right)\\ & \wedge \lcm(q_2,\ldots,q_n) = p\end{split}\right\}.\]

Again, we assume a JND of 18 cent between every two notes of a chord. Even though relative periodicity resembles harmonicity well qualitatively, Stolzenburg \cite{Stolzenburg2015Harmonyperceptionperiodicity} considers logarithmic periodicity as a computational model for harmonicity because of the Weber-Fechner law as discussed in Section \ref{sec:logperception}. We generalize JND to chords by determining a polyhedral neighborhood $N_c\subset \cC_n$ for each chord $c = [c_1,\ldots,c_n]\in \cC_n \setminus \cC_{n-1}$ in which there is no noticeable difference compared to $c$. Formally, we have for ${\rm JND} = 18$ cent
\begin{equation}\label{eq:jndnbhd}
	N_{c} :=\left\{[c_1+d_1,\ldots,c_n+d_n] \,\middle\vert\, \forall_{i,j=1,\ldots,n}\left(d_i\in [-{\rm JND},{\rm JND}] \wedge |d_i-d_j|\le {\rm JND}\right) \right\}.
\end{equation}

We could try generalizing periodicity to arbitrary chords using Table \ref{tab:periodicities}. Following Example 10 in \cite{Stolzenburg2015Harmonyperceptionperiodicity}, the first inversion of the diminished triad can be written as $[0,3,9]$, $[-3,0,6]$ and $[-9,-6,0]$. These representations have relative periodicities 15, 25 and 6 depending on which note in this triad is considered to be pitch 0 in Table \ref{tab:periodicities}. To illustrate the computation we compute frequency ratios $[5/6,1,7/5]$ from Table \ref{tab:periodicities} for $[-3,0,6]$ which translates to $[1,6/5,(6/5)\cdot (7/5)]$ for $[0,3,9]$ and results in an overall periodicity of $\lcm(1,5,5\cdot 5) = 5 \cdot 5 = 25$. In \cite{Stolzenburg2015Harmonyperceptionperiodicity} this problem of potentially having different periodicities for the same chord has been solved by computing the average of the three periodicities (both raw and logarithmic), i.e. raw $(15+25+6)/3 \approx 15.3$ and logarithmic $(\log_2(15)+\log_2(25)+\log_2(6))/3 \approx 3.7$. Even though this gives good empirical results, it seems to contradict the ``rational tuning'' principle, which uses the fractions with the smallest denominator approximating equal temperament within a certain error margin. For example, $19/16$, $13/11$ and $6/5$ all approximate the frequency ratios for the minor third, and $5$ is the relative periodicity. For the same reason, the relative periodicity for the the first inversion of the diminished triad should be 6 and not an average. In summary, the algorithm for finding the rational tuning provided by \cite{Stolzenburg2015Harmonyperceptionperiodicity} should be generalized to arbitrary chords rather than using the frequencies in Table \ref{tab:periodicities}.

Let us therefore modify the computation of periodicity slightly and not use the proposed smoothing from \cite{Stolzenburg2015Harmonyperceptionperiodicity}. We define for $c \in \cC^{[0,12]}_n$ 
\begin{equation}
	p(c):=\min\{p \mid N_c \cap \cC^p_n \neq \emptyset\}.\label{eq:chordperiodicity}
\end{equation}
Informally, we choose the best fit of periodicity for each chord within a JND for every two notes, rather than averaging over periodicities. 
Algorithm \ref{alg:periodicity} yields the periodicity of a chord with $n$ notes as a step function on $\cC^{[0,12]}_n$. We have used a resolution of 100 cent per semitone. For $n=4$ the array size of $\cC^{[0,12]}_n$ is therefore $1200^4 \approx 2\cdot 10^{12}$ which was our computational limit.\footnote{This can be implemented more efficiently, e.g. by only considering all $c\in \cC^p_n$ in a neighborhood of $RemainingChords$ and by reducing the resolution.}
\begin{algorithm}
	\caption{Determine periodicity step function $p:\cC_n^{[0,12]} \to \R$}\label{alg:periodicity}
	\begin{algorithmic}
		\Require $n \geq 1$ \Comment{$n$=number of chord tones}
		\State $q \gets 1$ \Comment{$q$=periodicity index}
		\State $RemainingChords \gets \cC^{[0,12]}_n$ \Comment{Consider all chords of the form $[0,c_2,\ldots,c_n]$}
		\While{$RemainingChords \neq \emptyset$} \Comment{While there are chords without periodicity}
		\ForAll{$c\in \cC^q_n$} \Comment{For all chords with periodicity $q$}
		\ForAll{$d\in N_c  \cap RemainingChords$} \Comment{For all new chords within JND}
		\State $p(d) \gets q$ \Comment{Set periodicity to $q$}
		\State $RemainingChords \gets (RemainingChords \setminus N_c)$ \Comment{Update new chords}
		\EndFor
		\EndFor
		\State {$q \gets q+1$} \Comment{Increase periodicity index by 1}
		\EndWhile
	\end{algorithmic}
\end{algorithm}

\begin{figure}[ht]
		\centering
	
	\includegraphics[width=0.9\textwidth,trim=50 53 20 70, clip]{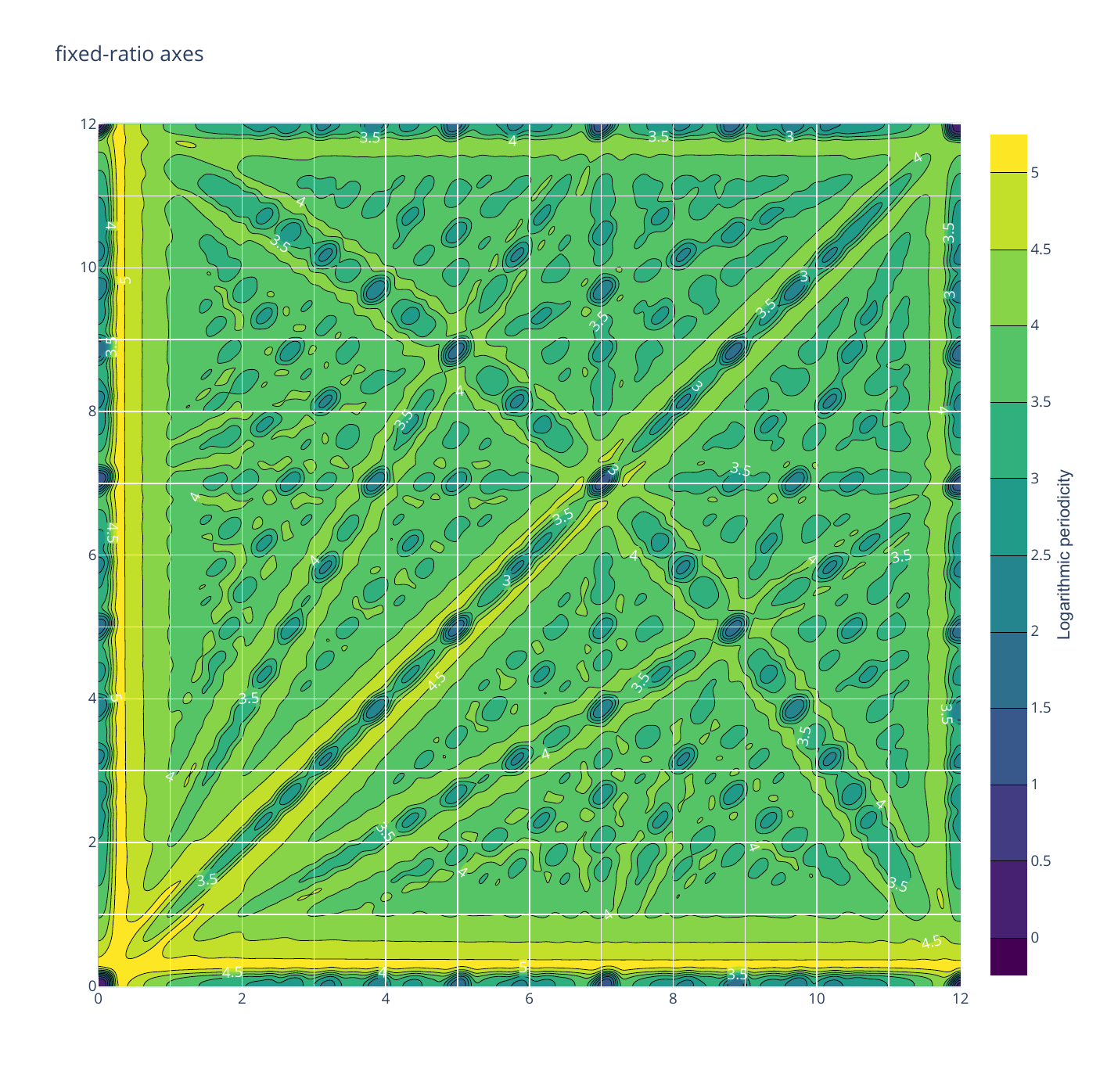}
	\caption{Contour plot for logarithmic periodicities of $[0,x,y]\in \cC_3^{[0,12]}$ using ${\rm JND} = 18$ cent smoothed by a Gaussian with standard deviation $\sigma = 6$ cent cent}\label{fig:gaussian_log_periodicity_triads_jnd18}
\end{figure}

This can be smoothed as discussed above using a Gaussian with standard deviation $6$ cent. Figure \ref{fig:gaussian_log_periodicity_triads_jnd18} visualizes the resulting logarithmic periodicity function $\log_2 p(c)$ for triads spanning at most one octave; we normalize a triad in continuous pitch space to be of the form $[0,x,y]$ with $x,y\in[0,12]$ and draw the graph as a contour plot in the $xy$-plane with the height $z$ given by the logarithmic periodicity. The intersection points of the grid lines correspond to chords in twelve-tone equal temperament, $\cC_2$ diagonally embeds into $\cC_3$, and the second inversion [0,5,9] of the major triad appears to be the most consonant chord consisting of 3 different tones. 
\captionsetup[figure]{skip=0pt}
\captionsetup[figure]{skip=10pt}

\section{Music perception}\label{sec:expectation}

Time adds a layer of complexity to sound perception by way of considering paths and sequences in $\cC$. We will focus on musical expectation as well as on tension and release. Our ability to anticipate future events is another vital aspect of human evolution. Perceptual expectation has been studied in cognitive neuroscience \cite{Bubic2009}, and it is also a fundamental part of music perception \cite{Huron2006,Pearce2012,Rohrmeier2013,Schmuckler1989,Seger2013,Wall2020}. Clearly, musical expectation depends on the listener, or, more precisely, on his brain and its musical training \cite{Bigand1996Perceptionmusicaltension}. It involves recognizing and predicting patterns both in sound and time set within a context. Among musicians this is also known as the concept of tension and release. Some aspects of its neuro-acoustic mechanism have been studied in the literature \cite{Chan2019ScienceHarmonyPsychophysical,TRAMO_2005,Lahdelma_2015}. While roughness plays a role in tuning, scales and the sound perception of chords, it is not audible in the psychoacoustic interaction between consecutive chords. We will demonstrate in this section why we should and how we can analyze a periodicity approach to tension and release using tools from differential geometry.

\subsection{Tension and release}\label{sec:resolutionsense}


Concepts related to the harmonic transitions are the circle of fifths, the Tonnetz model by Euler \cite{Euler2017,Euler:1774}, the tonal pitch space by Lerdahl \cite{Lerdahl.2001}, the tonal hierachy by Bharucha and Krumhansl \cite{Bharucha.1983, Krumhansl.2001}, the gauge-theoretic approach to tonal attraction \cite{Graben2017,Blutner2020}, and similar geometric structures describing harmonic relationships \cite{Cohn.1997,Callender2008,Tymoczko_2012}. Plus, the bass line clearly plays an important role in the perception of chord progressions \cite{Schwitzgebel2021,Hove2014}. There have been several studies about the relationship between these horizontal approaches and roughness, (vertical) harmonicity and as well as cultural familiarity \cite{Bigand_1999,Bigand1996Perceptionmusicaltension,Lahdelma2020Culturalfamiliaritymusical,Lahdelma_2016,Lahdelma_2015}. Both the models by Lehrdahl \cite{Lerdahl.2001} and by Bharucha and Krumhansl \cite{Bharucha.1983} capture and describe distances for harmonic motions, and can both be viewed as a metric space in the mathematical sense \cite{Randall2010}. We will not view perceptive distances in harmonic transitions as a metric, because the order of tones or chords matters. Instead, we will show how to make the relations between notes and chords depend on the context and the order, and thereby remedy the limitations of geometric models mentioned in \cite[119ff.]{Krumhansl.2001}. The intricate interplay between the voice leading distance presented in Section \ref{sec:pitchspace} and harmonic transition is important, even though harmonies and harmony theory are often discussed without considering the individual note movement. The cognitive mechanisms are related and interfere with each other to create the sensation of transitional harmony. This is hinted at in \cite{Tymoczko_2012,Callender2008}.

Let us discuss the potential factors that affect transitional harmonicity. On the one hand we expect the basic principles behind harmonicity to play a role. We have not found much empirical evidence, but we hypothesize that the neuronal network mechanisms behind many sensations and particularly between horizontal and vertical harmony should be the same, and Tramo et al. \cite[p. 96]{Tramo.2001} also suggest that the vertical and the horizontal dimension of harmony is related. Therefore phase locking or an even more fundamental physiological principle will be behind transitional harmony. On the other hand we expect pattern recognition and the universal ability of detecting minima in sensoric input to be important. Minimization is implicitly used in defining the voice leading distance presented in Section \ref{sec:pitchspace} given by the geodesic distance. 

For our purpose there are two fundamentally different kinds of expectations:
\begin{enumerate}
	\item If you listen to a piece of music, you can predict how it continues. You might be able to anticipate a few notes and chords depending on your training and background. Your anticipation will be based on tempo, meter, rhythm, melody, dynamics, form, chord progression which are time-dependent aspects of context. In the language of our geometric model: From a path in $\cC$ we want to anticipate its continuation. This will depend on its speed and its shape, including its direction, its curvature and other geometric aspects. You can compare a musical piece to a roller coaster ride, which you should construct or analyze using differential geometry. However, it will also depend on a second, time-independent kind of expectation.
	\item Assume you are listening to a single chord, and you have to predict which chord could follow. You might wonder, which context this is in, and this might partially be responsible for your expectations. As before, it is based on your training and background. However, there are physical reasons for your anticipations as well. This certainly has to do with consonance of chords and voice leading, but also with the order in which two different chords are played. We can view this expectation either as a psychoacoustic evaluation of difference vectors on $\cC$ or of ordered pairs of chords. We call this time-independent psychoacoustic quality for an ordered pair of chords the {\em resolve}. This time-independent quantity has been studied under different names in \cite{Geer1962connotationmusicalconsonance,Maher1976NeedResolutionRatings,Arthurs2017Perceptionisolatedchords,Lahdelma2020Culturalfamiliaritymusical}, but we would like to emphasize its dependence on its contextual reference by giving it this new name.
\end{enumerate}

A progression of notes and chords with or without additional bass notes can be described as a sequence of points in $\cC$, which can be viewed as a discretization of a path in $\cC$ parameterized by time. It can be approximated by a differentiable path. Either way, we can study differential geometric properties like speed, momentum, acceleration, and angular speed in order to analyze and understand chord progressions better. Furthermore, we can consider differential geometric properties of the path after applying suitable height functions. We hypothesize that the time-independent expectation can be deduced from the resolve by way of differential geometry. If we are at a point $p\in \cC$, we can quantify the resolve as a height function on $\cC$. In other words, the resolve is a function on $\cC \times \cC$.

Some interesting questions arise, which we do not attempt to answer here: Is the resolve the result of a priming with ordered pairs of chords based on cultural familiarity and training, or is it a multi-dimensional vector intrinsic to the starting chord or a local neighborhood of the starting chord, i.e. without the necessity of having ever heard the second chord? Is the training happening on the level of some basic neuronal mechanism for any chord progression or do we need all kinds of pairs of chords as training data? Do musicians and composers imagine the succeeding chord or do they sense the direction in which they have to move the notes? 

Research from \cite{Chan2019ScienceHarmonyPsychophysical} attempts to quantify the resolve. They call it transitional harmony and compute it via 
%
\begin{equation}
	\Delta\Delta \hat{t} := \frac{\Delta t_p-\Delta t_s}{T_{sub}}, \text{where } \Delta t = [k_it_i]_{max}-[k_it_i]_{min},\label{eq:transitionalharmony}
\end{equation}
where $[k_it_i]_{max}$ and $[k_it_i]_{min}$ are the largest and smallest multiples of the chord tone periods that (nearly) coincide with the chord periodicity $T_{sub}$ which we introduced in Section \ref{sec:chordperiodicity} and where the indices $s$ and $p$ correspond to the succeeding and nearest preceding chord, respectively. Even though the authors have found some strong correlations \cite[Table 3]{Chan2019ScienceHarmonyPsychophysical} supporting the validity of $\Delta\Delta \hat{t}$, we question the definition due to its strong dependence on small pitch changes: music perception should not change a lot by small pitch variations, but it does in the definition given by Equation \eqref{eq:transitionalharmony}. We hypothesize that the correlations found in \cite[Table 3]{Chan2019ScienceHarmonyPsychophysical} are due to the correlation between harmonicity and roughness for instruments with harmonic spectra. 


\subsection{Two-chord progressions starting with a tritone}\label{sec:tritone}

In order to motivate various approaches to the resolve we consider two-chord progressions of dyads within 12TET starting on a tritone $[F3,B3]$ where at least one note changes and each note does not move more than a semitone. Let us ignore the choice of octave in this section. There is a total of eight such chord movements.

The two parallel tritone movements considered on their own and out of context do not sound like they resolve anything, but adding the bass lines $C\sharp3\rightarrow F\sharp2$ or $G2\rightarrow C2$ yields the standard chord progression $II^7\rightarrow V^7$ as $[C\sharp3,E\sharp3,B3] \rightarrow [F\sharp2,E3,A\sharp3]$\footnote{The first note in this notation always corresponds to the bass note.} and $[G2,F3,B3]\rightarrow[C2,E3,B\flat3]$, respectively, where the tonalities are clearly very far away from each other. On the one hand this simple example confirms the well-known assumption that chords should always be viewed in a context, but on the other hand it represents the charm behind the technique of modulation in music. It is therefore nevertheless necessary to consider chord progressions without a given tonality or context. It just leaves chord progressions ambiguous, and probabilistic methods can be employed.

The strongest resolution from the perspective of periodicities or ratios should clearly the progression to the perfect fifth $[F3,B3] \rightarrow [E3,B3]$ or $[F3,B3] \rightarrow [F3,C3]$. However, it does not sound like a good way of resolving the tritone. If we think of the notes as attracting or repelling magnets then both the notes should move in opposing directions in order to resolve the dissonance, which we will consider in the next paragraph. However, we can again add bass lines to make the first progression sound like the jazz resolution to the major seventh chord $V^7\rightarrow I^{\Delta}$ given by $[G2,F3,B3] \rightarrow [C2,E3,B3]$ and the second progression like the resolution $V^{\text{dim}7}\to I^\Delta$ or $V^7\to I^{\text{sus}4} \to I$ partially represented by $[A\flat2,F3,C\flat4] \rightarrow [D\flat2,F3,C4]$ and $[G2,F3,B3] \rightarrow [C2,F3,C4]\rightarrow [C2,E3,C4]$, respectively.

The chord progression into a perfect fourth $[F3,B3] \rightarrow [F3,Bb3]$ or $[F3,B3] \rightarrow [F\sharp3,B3]$ also does not sound like a good way of resolving the tritone. Again, we can put them in a suitable context by adding bass lines. The first progression sounds like the jazz resolution to the major seventh chord $V^7\rightarrow I^{\Delta}$ partially given by $[D\flat,F,B] \rightarrow [G\flat2,F3,Bb3]$ and the second progression like the resolution $V^{\text{dim}7}\to I^\Delta$ or $V^7\to I^{\text{sus}4} \to I$ partially represented by $[D3,F3,B3] \rightarrow [G2,F\sharp3,B3]$ and $[G2,F3,B3] \rightarrow [C2,F3,C4]\rightarrow [C2,E3,C4]$, respectively.

The best sounding dyad progression is the tritone $[F3,B3]$ resolving into the major third $[G\flat3,B\flat3]$ or the minor sixth $[E3,C4]$. Even though these progressions already sound like resolutions, it helps to view them in a context and a tonality in order to relate them to music theory. Possibly, our brain has already been primed for possible tonalities, and some tonalities are more probable than others. Clearly, the corresponding chord progressions are $V^7\to I$ partially represented by $[D\flat3,F3,C\flat4] \to [G\flat2,G\flat3,B\flat3]$ and $[G2,F3,B3] \to [C2,E3,C4]$. Notice that the tonality is already determined by the progression of dyads, the bass line only emphasizes the tonal center. An insightful work by Tom Sutcliffe \cite{Sutcliffe.2011} picks up on the gap in the literature of failing to explain why voice leading in combination with root progressions is used in tonal pieces. 



Let us describe a few possible approaches to transitional harmony between two chords, consider the differential geometry and revisit the above example.

\subsection{Transitive periodicity from the first to the second chord}\label{sec:transitive_periodicity}

Musical structures like rhythmic patterns and periodicity cause phase locking \cite{Nederlanden2020}. We therefore assume that the brain relates two chords $c_1$ and $c_2$ of a chord progression $c_1\rightarrow c_2$ through the working memory based on phase locking as described in Section \ref{sec:phaselocking}. On the one hand this seems compatible with the strong preference to descending fifths and ascending fourths over descending fourths and ascending fifths. On the other hand, if $c_1$ and $c_2$ only consist of one note each, a low periodicity of $c_1$ with respect to $c_2$ is desirable, because the neuronal firing is synchronized. This seems to be incompatible with voice leading at first, but as soon as you consider small chord movements with respect to voice leading, it is possible to move a short distance while being close with respect to phase synchronization. Therefore, we introduce a transitive periodicity analogously to the periodicity definition given in Equation \eqref{eq:chordperiodicity}. 

\begin{figure}[ht]
	\centering
	\def\svgwidth{0.3\textwidth}
\executeiffilenewer{Chord-Progression.svg}{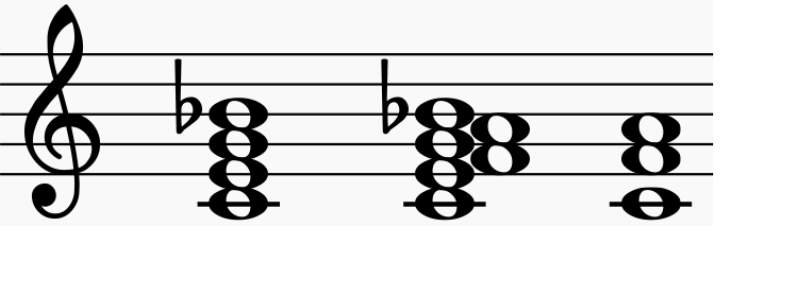}%
{inkscape -z -D --file=Chord-Progression.svg %
--export-pdf=Chord-Progression.pdf --export-latex}%
\begingroup%
  \makeatletter%
  \providecommand\color[2][]{%
    \errmessage{(Inkscape) Color is used for the text in Inkscape, but the package 'color.sty' is not loaded}%
    \renewcommand\color[2][]{}%
  }%
  \providecommand\transparent[1]{%
    \errmessage{(Inkscape) Transparency is used (non-zero) for the text in Inkscape, but the package 'transparent.sty' is not loaded}%
    \renewcommand\transparent[1]{}%
  }%
  \providecommand\rotatebox[2]{#2}%
  \newcommand*\fsize{\dimexpr\f@size pt\relax}%
  \newcommand*\lineheight[1]{\fontsize{\fsize}{#1\fsize}\selectfont}%
  \ifx\svgwidth\undefined%
    \setlength{\unitlength}{388.47774952bp}%
    \ifx\svgscale\undefined%
      \relax%
    \else%
      \setlength{\unitlength}{\unitlength * \real{\svgscale}}%
    \fi%
  \else%
    \setlength{\unitlength}{\svgwidth}%
  \fi%
  \global\let\svgwidth\undefined%
  \global\let\svgscale\undefined%
  \makeatother%
  \begin{picture}(1,0.37312816)%
    \lineheight{1}%
    \setlength\tabcolsep{0pt}%
    \put(0,0){\includegraphics[width=\unitlength,page=1]{Chord-Progression.pdf}}%
    \put(0.24467413,0.01562405){\color[rgb]{0,0,0}\makebox(0,0)[lt]{\lineheight{1.25}\smash{\begin{tabular}[t]{l}$c_1$\end{tabular}}}}%
    \put(0.77157057,0.01753918){\color[rgb]{0,0,0}\makebox(0,0)[lt]{\lineheight{1.25}\smash{\begin{tabular}[t]{l}$c_2$\end{tabular}}}}%
    \put(0.47087497,0.0147812){\color[rgb]{0,0,0}\makebox(0,0)[lt]{\lineheight{1.25}\smash{\begin{tabular}[t]{l}$[c_1,c_2]$\end{tabular}}}}%
  \end{picture}%
\endgroup%

	\caption{A chord progression $c_1\rightarrow c_2$ with the combined chord $[c_2,c_1]$}
	\label{fig:chordprogression}
\end{figure}

The transitive periodicity from $c_1$ to $c_2$ is the number of periods of $c_2$ necessary to match up with a period multiple of $c_1$, where the periods for $c_1$ and $c_2$ are each due to phase locking. Formally, transitive periodicity from $c_1$ to $c_2$ is the periodicity of $[c_1,c_2]$ relative to $c_2$, where $[c_1,c_2]$ is the (set-theoretic) union of $c_1$ and $c_2$: Figure \ref{fig:chordprogression} shows 
$c_1=[0,4,7,10]$, $c_2=[0,5,9]$ and the combined chord $[c_1,c_2]=[0,4,5,7,9,10]$. 

This corresponds to computing $p([c_1,c_2])/p(c_2)$, but due to the JND from Section \ref{eq:chordperiodicity} the smoothed periodicities are not the correct quantities to be used for computing transitive periodicity. We need to view $c_1$, $c_2$ and $[c_1,c_2]$ in the context of their related periodicities before smoothing. Due to technical difficulties we need to work with $\R^n$ rather than its quotient $\cC^n$.  
In analogy to Section \ref{sec:chordperiodicity} we define
\begin{equation}\cC^{p}_{m,n} := \left\{\left(0,12\cdot\log_2\left(\tfrac{p_2}{q_2}\right),\ldots,12\cdot\log_2\left(\tfrac{p_{m+n}}{q_{m+n}}\right)\right) \,\middle\vert\,\\  \begin{split}&\forall_{i} \left(1 \le \tfrac{p_i}{q_i}\le 2 \wedge \gcd(p_i,q_i) = 1\right)\\ & \wedge \frac{\lcm(q_1,\ldots,q_{m+n})}{\lcm(q_2,\ldots,q_n)} = p\end{split}\right\}\label{eq:transitive-pSpace}\end{equation}
and for pitch tuples $t \in\R^n$
\begin{equation}\label{eq:transitiv_jndnbhd}
	N_{t} :=\left\{t+(d_1,\ldots,d_n) \,\middle\vert\, \forall_{i,j=1,\ldots,n}\left(d_i\in [-{\rm JND},{\rm JND}] \wedge |d_i-d_j|\le {\rm JND}\right) \right\}.
\end{equation}
Informally, $\cC^p_{m,n}$ contains the $(m+n)$--tuples $(t_1,t_2)$ with representatives $t_1 \in \R^m$ und $t_2\in \R^n$ of $c_1 \in \cC_m$ and $c_2 \in \cC_n^{[0,12]}$ so that the periodicity of $[c_1,c_2]$ relative to $c_2$ is $p$.\footnote{In order to improve readability the first $n$ entries in the elements of $\cC^p_{m,n}$ correspond to $c_2$.}
In order to incorporate approximations within a JND we define transitional periodicity $p:\cC_m \times \cC_n^{[0,12]} \to \R, (c_1,c_2) \mapsto p(c_1 \to c_2)$ via
\[
p(c_1 \to c_2):=\min\left\{p \,\middle\vert\,\\  \begin{split} &N_t \cap \cC^p_{m,n} \neq \emptyset \text{, where } t \in \R^{m+n},\\ &c_1 = [t_1,\ldots,t_m], c_2 = [t_{m+1},\ldots,t_{m+n}]\end{split}\right\}.
\]

Just like in the case of periodicity in Section \ref{sec:chordperiodicity} we can use the logarithmic transitive periodicity. In order to extend $p$ to $\cC_m \times \cC_n$, it will be necessary to shift a chord $c_1$ and $c_2$ so that the lowest note of $c_2$ is $0$. For $c=[p_1,\ldots,p_n]$ let $s_p(c) := [p_1-p,p_2-p,\ldots,p_n-p]$ and $\min(c):=\min\{p_1,\ldots,p_n\}$. Then $p(c_1\to c_2)$ will be redefined as $p(s_{\min (c_2)}(c_1),s_{\min (c_2)}(c_2))$.

\begin{figure}[ht]
	\centering
	\includegraphics[width=0.47\textwidth,trim=50 53 20 70, clip]{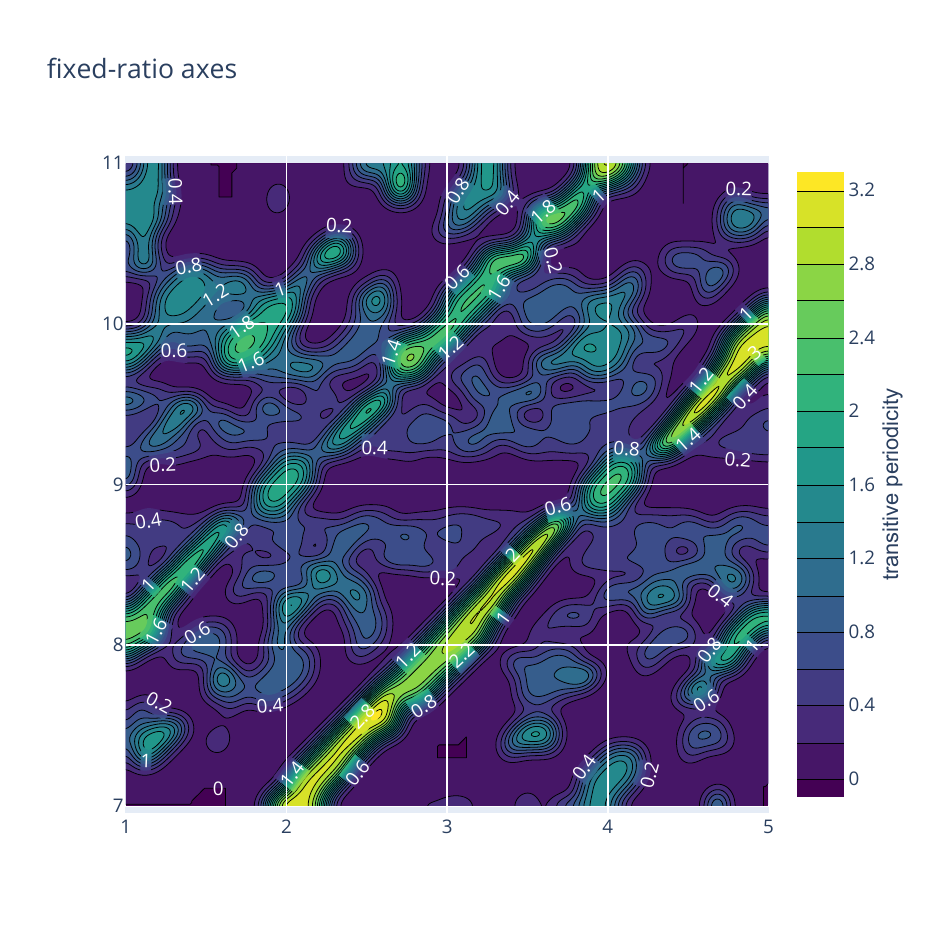}
	\includegraphics[width=0.47\textwidth,trim=50 53 20 70, clip]{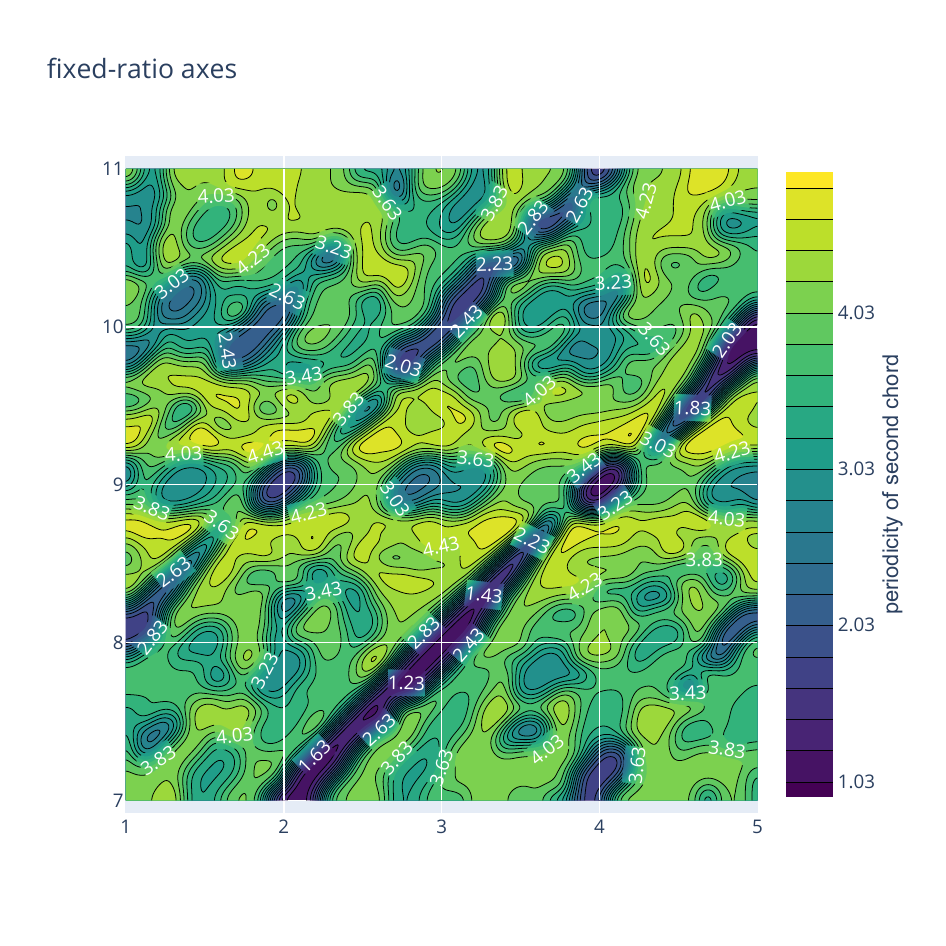}
	\caption{Contour plot for logarithmic transitive periodicities $p(c_1\to c_2)$ (left) with the corresponding logarithmic periodicities $p_{c_1}(c_2)$ of $c_2\in\cC_2$ (right) starting with the tritone $c_1=[3,9]$ using ${\rm JND} = 18$ cent smoothed by a Gaussian with standard deviation $\sigma = 6$ cent}\label{fig:new_periodicity_dyads}
\end{figure}

Let us revisit the example in Section \ref{sec:tritone}. The local neighborhood of the transitional periodicity $c_1\to c_2$ starting with the tritone $c_1=[3,9]$ with the corresponding periodicities of $c_2$ is shown in Figure \ref{fig:new_periodicity_dyads}.

\captionsetup[figure]{skip=0pt}
\captionsetup[figure]{skip=10pt}

Notice that while the periodicities of the perfect fourth and fifth are smallest, the transitive periodicities resolving to the perfect fourth and fifth are bigger than the transitive periodicities resolving to other chords. Even if small transitive periodicities play a role in chord resolutions, they do not fully explain them.
The algorithm for determining the transitive periodicity step function as in Equation \eqref{eq:chordperiodicity} is shown in Algorithm \ref{alg:transitive_periodicity}, where we define for $c =[c_1,\ldots,c_m]\in \cC_m$ \[N_c(\epsilon,n) :=\left\{(t_1+d_1,\ldots,t_n+d_n) \,\middle\vert\, t \in \R^n, [t]=c,\forall_{i=1,\ldots,n}\left(d_i\in [-\epsilon,\epsilon] \right) \right\} \subset \cC_n,\] since we want $c_1$ and $c_2$ to be close with respect to voice leading. In the algorithm we need the projection to the last $n$ coordinates 
$l_n(t_1,\ldots,t_{m+n}):=(t_{m+1},\ldots,t_{m+n})$.
\begin{algorithm}
	\caption{Determine transitive periodicity step function $p:\cC \times \cC \to \R$\label{alg:transitive_periodicity}}
	\begin{algorithmic}
		\Require $m,n \geq 1$ \Comment{$m,n$=number of tones in $c_1,c_2\in \cC$}
		\Require $c_1 \in \cC_m$ \Comment{$c_1$=chord with $m$ tones}
		\Require $\epsilon \ge 0$ \Comment{$\epsilon$=maximal distance between the tones of $c_1 \in \cC_m$ and $c_2\in \cC_n$}
		\State $q \gets 1$ \Comment{$q$=periodicity index}
		\State $RemainingChords \gets N_{c_1}(\epsilon,n)$ \Comment{Consider all chords with $n$ notes close enough to $c_1$}
		\While{$RemainingChords \neq \emptyset$} \Comment{While there are chords without periodicity}
		\ForAll{$s\in \cC^q_{m,n}$} \Comment{For all chords with relative periodicity $q$}
		\ForAll{$t \in N_{l_n(s)}  \cap RemainingChords$} \Comment{For all new chords within JND}
		\State $p(c_1,[t]) \gets q$ \Comment{Set periodicity to $q$}
		\State $RemainingChords \gets (RemainingChords \setminus N_{l_n(s)})$ \Comment{Update new chords}
		\EndFor
		\EndFor
		\State {$q \gets q+1$} \Comment{Increase periodicity index by 1}
		\EndWhile
	\end{algorithmic}
\end{algorithm}

We hypothesize that transitive periodicity will play a role in combination with the usual periodicity. A good chord resolution will have a low periodicity for the second chord $c_2$ as well as a low transitive periodicity $c_1 \to c_2$.

The combined chord offers more possibilities for transitive quantities that can be studied in the context of music perception. For example, we can consider the periodicity $p_{c_1}([c_1,c_2])$ of $[c_1,c_2]$ relative to $c_1$ for the chord progression $c_1 \to  c_2$.

\subsection{Directional derivative of periodicity}

Melodies and chord progressions not only have a sense of direction in $\cC$, but the rate of change in psychoacoustic quantities with respect to these directions should play a role in music perception. Assuming that periodicity is a good measure for consonance, chord resolutions should decrease periodicity while traveling only a small distance with respect to voice leading. If this perceptive quality is truly local then the infinitesimal change in periodicity can be formulated using directional derivatives in chord space.

Due to the soft computing skills of our brain, periodicity on $\cC$ can be considered to be continuously differentiable. Let $v$ be a tangent vector of $T_c\cC$ at a chord $c\in \cC$. Let $\rho:[-\epsilon,\epsilon]\to \cC$ be a continously differentiable path with $\dot \rho|_{t=0} = v$. Then the directional derivative of the periodicity in the direction $v \in T_c \cC$ at $c\in \cC$ is given by 
\[
D_v p:=\left.\tfrac{d}{dt}\right|_{t=0} (p \circ \rho).\]
The more negative $D_v p$ is, the stronger its resolution in the direction $v$ is perceived. It is not enough that the final chord is more consonant. For example, either resolution from a tritone to a perfect fifth in 12TET does not sound as good as the one to the minor sixth or the major third as we have discussed in Section \ref{sec:tritone}.

On the other hand, the parallel movement of chords does not change periodicity. This implies that $D_v p$ vanishes for $v=[1,1]$. Clearly, maximal infinitesimal change in periodicity (either negative or positive) is provided for $v=[-1,1]$ and $v=[1,-1]$. If Figure \ref{fig:periodicity-dyads-smoothed} is correct, however, then only a very small repelling movement will reduce periodicity which does not yield harmonic relationships of chords. The progressions $[3,9] \to [2.75,9.25]$ and $[3,9] \to [3.25,8.75]$ will increase periodicity. The progressions to the perfect fourth and fifth $[3,9] \to [2.5,9.5]$ and $[3,9] \to [3.5,8.5]$ will certainly decrease periodicity. Interestingly enough, the only feasible progressions $[3,9] \to [2,10]$ and $[3,9] \to [4,8]$ in 12TET do not reduce periodicity by much or at all. Still, they are the best resolutions available in 12TET.

While the directional derivative presents an interesting approach it needs to be viewed in combination with other aspects of transitional harmony as part of a Pareto optimal solution. Possibly, the periodicity function shown in Figure \ref{fig:periodicity-dyads-smoothed} needs to be corrected as well. However, it shows that, in the case of the tritone, the chord progressions with the maximal effect on periodicity will move each note simultaneously inwards or outwards. Furthermore, it suggests that quarter tone movements will be the best chord resolution when considering periodicity only. This hypothesis is confirmed by the author's perception.

\section{Results and Discussion}


Music is considered as something real and vital for humankind, but no attempt on a holistic model for music perception has yet been attempted. Clearly, humans do not perceive musical sounds as the complicated audio waves they are or as the way it is presented in music notation but as something simple and often beautiful. While music theory formalizes the music we perceive, music psychology carries out empirical studies about specific perceptive aspects. We envision, that it should also be possible to deduce music theory from music psychology with the correct holistic model for music perception. Using psychoacoustic results and facts from music theory it should be possible to reverse engineer this model. With this in mind, we have introduced mathematical structures that allow for rigorous quantitative studies of music perception based on the mechanics described by physical or neuronal models. We laid an emphasis on a rigorous approach that is not more complicated than absolutely necessary and which can be extended when needed.

We revisit ideas by Tymoczko \cite{Tymoczko.2006,Tymoczko.2011} to prove that the space of chords $\cC$ is a metric space and a Whitney stratified space with a Riemannian structure. The geometry of $\cC$ is not much more than Euclidean space itself. However, it allows us to apply calculus across different strata of $\cC$. Furthermore, the Riemannian metric on $\cC$ allows us to consider the geodesic distance across different strata which yields a voice leading distance satisfying the triangle inequality. The geodesic approach is surprisingly simple and natural considering the common desire that distance functions satisfy certain conditions and in view of more elaborate attempts regarding voice leading distances \cite{Callender2008,Milne2016,Genuys2019}. The space $\cC$ only contains the objects for music production, but not any information about music perception. 

Psychoacoustic quantities can be viewed, computed and analyzed as height functions on $\cC$. In particular, we have modified the periodicity approach to consonance by Stolzenburg \cite{Stolzenburg2015Harmonyperceptionperiodicity} in order to present a definition of periodicity for arbitrary chords. Roughness is another way of interpreting dissonance. Height functions themselves are static. Music is a dynamic process, so it might be necessary to consider the change in height function as a dynamical system in order to deduce properties of music. All of the psychoacoustic functions can be assumed to be differentiable which enables us to use tools from differential geometry to study them by considering gradient vectors and directional derivatives.

The height function for periodicity led to two possible approaches for transitive harmonicity. In particular, we showed how to use the differential structure of the periodicity graph on $\cC$ to study geometric properties of paths in $\cC$ and their respective lifts to the graph of psychoacoustic functions on $\cC$. We implicitly assume that the geodesic distance agrees with the psychoacoustic reality. This needs to be verified empirically. Although we do not expect our two approaches for transitive harmonicity to be valid, we expect that other approaches to music perception can be analyzed using the differential geometric framework. Clearly, music works, because we look at a discrete subset of chords with certain properties. It will be interesting to see which tools are the correct ones for discretizing the differential-geometric model.

The differential-geometric structure invites studies that falsify or confirm psychoacoustic models for music. Ultimately, this approach can close the gap between music theory and music psychology. Even though the mechanisms discussed here stem from Western music, they are founded on more general physical and neuronal principles, which are in theory applicable to music from other cultures or sounds with inharmonic spectra. Furthermore, it will be interesting to study, generalize and extend the mathematical structures themselves and to incorporate statistical aspects of music perception in the model.

\section{Conclusion}

For the purpose of analyzing music perception, we have described useful geometric structures for the space of chords $\cC$. We have rigorously proven properties that are desirable from a mathematical as well as from a music perception point of view. In particular, chords with a different number of notes can be viewed as strata of $\cC$. The Riemannian metric on each stratum allows to define a geodesic distance on $\cC$, which makes it into a metric space. The metric is a natural choice for determining efficient voice leading. The Riemannian metric also allows to study shapes in the context of music perception. This enables music psychologists and music theorists to use tools from differential geometry in order to study music perception.

\bibliography{references2}

\begin{thebibliography}{137}
\providecommand{\natexlab}[1]{#1}
\providecommand{\url}[1]{\texttt{#1}}
\expandafter\ifx\csname urlstyle\endcsname\relax
  \providecommand{\doi}[1]{doi: #1}\else
  \providecommand{\doi}{doi: \begingroup \urlstyle{rm}\Url}\fi

\bibitem[Abdi(2010)]{Abdi2010}
H.~Abdi.
\newblock Signal detection theory.
\newblock In \emph{International Encyclopedia of Education}, pages 407--410.
  Elsevier, 2010.
\newblock \doi{10.1016/b978-0-08-044894-7.01364-6}.

\bibitem[Akkoc(2002)]{Akkoc2002}
C.~Akkoc.
\newblock Non-deterministic scales used in traditional turkish music.
\newblock \emph{Journal of New Music Research}, 31\penalty0 (4):\penalty0
  285--293, dec 2002.
\newblock \doi{10.1076/jnmr.31.4.285.14169}.

\bibitem[Alekseevsky et~al.(2003)Alekseevsky, Kriegl, Losik, and
  Michor]{Alekseevsky2003RiemannianGeometryOrbit}
D.~Alekseevsky, A.~Kriegl, M.~Losik, and P.~W. Michor.
\newblock The riemannian geometry of orbit spaces. the metric, geodesics, and
  integrable systems.
\newblock \emph{Publ. Math. Debrecen}, 6\penalty0 (3--4):\penalty0 247--276,
  2003.
\newblock URL \url{https://publi.math.unideb.hu/load_doc.php?p=841&t=pap}.

\bibitem[Armitage et~al.(2021)Armitage, Lahdelma, and Eerola]{Armitage2021}
J.~Armitage, I.~Lahdelma, and T.~Eerola.
\newblock Automatic responses to musical intervals: Contrasts in acoustic
  roughness predict affective priming in western listeners.
\newblock \emph{The Journal of the Acoustical Society of America}, 150\penalty0
  (1):\penalty0 551--560, jul 2021.
\newblock \doi{10.1121/10.0005623}.

\bibitem[Arthurs et~al.(2017)Arthurs, Beeston, and
  Timmers]{Arthurs2017Perceptionisolatedchords}
Y.~Arthurs, A.~V. Beeston, and R.~Timmers.
\newblock Perception of isolated chords: Examining frequency of occurrence,
  instrumental timbre, acoustic descriptors and musical training.
\newblock \emph{Psychology of Music}, 46\penalty0 (5):\penalty0 662--681, aug
  2017.
\newblock \doi{10.1177/0305735617720834}.

\bibitem[Bailes et~al.(2015)Bailes, Dean, and Broughton]{Bailes2015}
F.~Bailes, R.~T. Dean, and M.~C. Broughton.
\newblock How different are our perceptions of equal-tempered and microtonal
  intervals? a behavioural and {EEG} survey.
\newblock \emph{{PLOS} {ONE}}, 10\penalty0 (8):\penalty0 e0135082, aug 2015.
\newblock \doi{10.1371/journal.pone.0135082}.

\bibitem[Balkwill and Thompson(1999)]{Balkwill.1999}
L.-L. Balkwill and W.~F. Thompson.
\newblock {A Cross-Cultural Investigation of the Perception of Emotion in
  Music: Psychophysical and Cultural Cues}.
\newblock \emph{{Music Perception: An Interdisciplinary Journal}}, 17\penalty0
  (1):\penalty0 43--64, 1999.
\newblock ISSN 0730-7829.
\newblock \doi{10.2307/40285811}.
\newblock URL \url{https://mp.ucpress.edu/content/17/1/43}.

\bibitem[Bausenhart et~al.(2018)Bausenhart, Luca, and Ulrich]{Bausenhart2018}
K.~M. Bausenhart, M.~D. Luca, and R.~Ulrich.
\newblock Assessing duration discrimination: Psychophysical methods and
  psychometric function analysis.
\newblock In \emph{Timing and Time Perception: Procedures, Measures, {\&}
  Applications}, pages 52--78. {BRILL}, mar 2018.
\newblock \doi{10.1163/9789004280205_004}.

\bibitem[Becker(2019)]{Becker2019}
J.~Becker.
\newblock \emph{Traditional Music in Modern Java}.
\newblock University of Hawaii Press, mar 2019.
\newblock \doi{10.2307/j.ctv9zcjt8}.

\bibitem[beim Graben and Blutner(2017)]{Graben2017}
P.~beim Graben and R.~Blutner.
\newblock Toward a gauge theory of musical forces.
\newblock In \emph{Quantum Interaction}, pages 99--111. Springer International
  Publishing, 2017.
\newblock \doi{10.1007/978-3-319-52289-0_8}.

\bibitem[Bettiol et~al.(2018)Bettiol, Derdzinski, and Piccione]{Bettiol2018}
R.~G. Bettiol, A.~Derdzinski, and P.~Piccione.
\newblock Teichmüller theory and collapse of flat manifolds.
\newblock \emph{Annali di Matematica Pura ed Applicata (1923 -)}, 197\penalty0
  (4):\penalty0 1247--1268, jan 2018.
\newblock \doi{10.1007/s10231-017-0723-7}.

\bibitem[Bharucha and Krumhansl(1983)]{Bharucha.1983}
J.~Bharucha and C.~L. Krumhansl.
\newblock {The representation of harmonic structure in music: Hierarchies of
  stability as a function of context}.
\newblock \emph{{Cognition}}, 13\penalty0 (1):\penalty0 63--102, 1983.
\newblock ISSN 0010-0277.
\newblock \doi{10.1016/0010-0277(83)90003-3}.
\newblock URL
  \url{http://www.sciencedirect.com/science/article/pii/0010027783900033}.

\bibitem[Bidelman and Krishnan(2009)]{Bidelman2009NeuralCorrelatesConsonance}
G.~M. Bidelman and A.~Krishnan.
\newblock Neural correlates of consonance, dissonance, and the hierarchy of
  musical pitch in the human brainstem.
\newblock \emph{Journal of Neuroscience}, 29\penalty0 (42):\penalty0
  13165--13171, oct 2009.
\newblock \doi{10.1523/jneurosci.3900-09.2009}.

\bibitem[Bigand and Parncutt(1999)]{Bigand_1999}
E.~Bigand and R.~Parncutt.
\newblock Perceiving musical tension in long chord sequences.
\newblock \emph{Psychological Research}, 62\penalty0 (4):\penalty0 237--254,
  oct 1999.
\newblock \doi{10.1007/s004260050053}.

\bibitem[Bigand et~al.(1996)Bigand, Parncutt, and
  Lerdahl]{Bigand1996Perceptionmusicaltension}
E.~Bigand, R.~Parncutt, and F.~Lerdahl.
\newblock Perception of musical tension in short chord sequences: The influence
  of harmonic function, sensory dissonance, horizontal motion, and musical
  training.
\newblock \emph{Perception {\&} Psychophysics}, 58\penalty0 (1):\penalty0
  125--141, jan 1996.
\newblock \doi{10.3758/bf03205482}.

\bibitem[Blutner and beim Graben(2020)]{Blutner2020}
R.~Blutner and P.~beim Graben.
\newblock Gauge models of musical forces.
\newblock \emph{Journal of Mathematics and Music}, 15\penalty0 (1):\penalty0
  17--36, feb 2020.
\newblock \doi{10.1080/17459737.2020.1716404}.

\bibitem[Borzellino(1992)]{Borzellino1992}
J.~E. Borzellino.
\newblock \emph{Riemannian Geometry of Orbifolds}.
\newblock PhD thesis, University of California, Los Angelos, 1 1992.

\bibitem[Boulos(2021)]{Boulos2021}
I.~Boulos.
\newblock Inside arabic music: Arabic maqam performance and theory in the 20th
  century. by johnny farraj and sami abu shumays.
\newblock \emph{Music and Letters}, 102\penalty0 (1):\penalty0 171--172, feb
  2021.
\newblock \doi{10.1093/ml/gcab018}.

\bibitem[Bridges(2008)]{Bridges2008}
B.~Bridges.
\newblock {Can Harmony be Non-Linear? a response to some of Glenn Branca’s
  ‘25 Questions’}.
\newblock In \emph{{Society for Musicology in Ireland annual conference}},
  Waterford Institute of Technology, 2008.

\bibitem[Bubic et~al.(2009)Bubic, von Cramon, Jacobsen, Schröger, and
  Schubotz]{Bubic2009}
A.~Bubic, D.~Y. von Cramon, T.~Jacobsen, E.~Schröger, and R.~I. Schubotz.
\newblock Violation of expectation: Neural correlates reflect bases of
  prediction.
\newblock \emph{Journal of Cognitive Neuroscience}, 21\penalty0 (1):\penalty0
  155--168, jan 2009.
\newblock \doi{10.1162/jocn.2009.21013}.

\bibitem[Burrows(1997)]{Burrows1997}
D.~Burrows.
\newblock A dynamical systems perspective on music.
\newblock \emph{Journal of Musicology}, 15\penalty0 (4):\penalty0 529--545,
  1997.
\newblock \doi{10.2307/764006}.

\bibitem[Burton(2015)]{Burton.2015}
R.~L. Burton.
\newblock {The Elements of Music: What Are They, and Who Cares?}
\newblock In J.~Rosevear and S.~Harding, editors, \emph{{ASME XXth National
  Conference Proceedings}}, Parkville, Victoria, 2015. {The Australian Society
  for Music Education Inc}.

\bibitem[Callender et~al.(2008)Callender, Quinn, and Tymoczko]{Callender2008}
C.~Callender, I.~Quinn, and D.~Tymoczko.
\newblock Generalized voice-leading spaces.
\newblock \emph{Science}, 320\penalty0 (5874):\penalty0 346--348, 2008.
\newblock ISSN 0036-8075.
\newblock \doi{10.1126/science.1153021}.
\newblock URL \url{https://science.sciencemag.org/content/320/5874/346}.

\bibitem[Cariani and Delgutte(1996)]{Cariani1996}
P.~A. Cariani and B.~Delgutte.
\newblock Neural correlates of the pitch of complex tones. i. pitch and pitch
  salience.
\newblock \emph{Journal of Neurophysiology}, 76\penalty0 (3):\penalty0
  1698--1716, sep 1996.
\newblock \doi{10.1152/jn.1996.76.3.1698}.

\bibitem[Chan et~al.(2019)Chan, Dong, and
  Li]{Chan2019ScienceHarmonyPsychophysical}
P.~Y. Chan, M.~Dong, and H.~Li.
\newblock The science of harmony: A psychophysical basis for perceptual
  tensions and resolutions in music.
\newblock \emph{Research}, 2019:\penalty0 1--22, sep 2019.
\newblock \doi{10.34133/2019/2369041}.

\bibitem[Cohen(1984)]{Cohen1984}
E.~A. Cohen.
\newblock Some effects of inharmonic partials on interval perception.
\newblock \emph{Music Perception}, 1\penalty0 (3):\penalty0 323--349, 1984.
\newblock \doi{10.2307/40285264}.

\bibitem[Cohn(1997)]{Cohn.1997}
R.~Cohn.
\newblock {Neo-Riemannian Operations, Parsimonious Trichords, and their
  'Tonnetz' Representations}.
\newblock \emph{{Journal of Music Theory}}, 42\penalty0 (2):\penalty0 1--66,
  1997.

\bibitem[Collins et~al.(2014)Collins, Tillmann, Barrett, Delb{\'{e}}, and
  Janata]{Collins2014}
T.~Collins, B.~Tillmann, F.~S. Barrett, C.~Delb{\'{e}}, and P.~Janata.
\newblock A combined model of sensory and cognitive representations underlying
  tonal expectations in music: From audio signals to behavior.
\newblock \emph{Psychological Review}, 121\penalty0 (1):\penalty0 33--65, 2014.
\newblock \doi{10.1037/a0034695}.

\bibitem[Cook and Fujisawa(2006)]{Cook2006PsychophysicsHarmonyPerception}
N.~D. Cook and T.~X. Fujisawa.
\newblock The psychophysics of harmony perception: Harmony is a three-tone
  phenomenon.
\newblock \emph{Empirical Musicology Review}, 1\penalty0 (2):\penalty0
  106--126, 2006.
\newblock \doi{10.18061/1811/24080}.

\bibitem[del Pozo and G{\'{o}}mez-Mart{\'{\i}}n(2022{\natexlab{a}})]{Pozo2022}
I.~del Pozo and F.~G{\'{o}}mez-Mart{\'{\i}}n.
\newblock A mathematical model of tonal function (i): Voice leadings.
\newblock In \emph{Mathematics and Computation in Music}, pages 218--230.
  Springer International Publishing, 2022{\natexlab{a}}.
\newblock \doi{10.1007/978-3-031-07015-0_18}.

\bibitem[del Pozo and G{\'{o}}mez-Mart{\'{\i}}n(2022{\natexlab{b}})]{Pozo2022a}
I.~del Pozo and F.~G{\'{o}}mez-Mart{\'{\i}}n.
\newblock A mathematical model of tonal function ({II}): Modulation.
\newblock In \emph{Mathematics and Computation in Music}, pages 231--239.
  Springer International Publishing, 2022{\natexlab{b}}.
\newblock \doi{10.1007/978-3-031-07015-0_19}.

\bibitem[Demorest et~al.(2009)Demorest, Morrison, Stambaugh, Beken, Richards,
  and Johnson]{Demorest.2009}
S.~M. Demorest, S.~J. Morrison, L.~A. Stambaugh, M.~Beken, T.~L. Richards, and
  C.~Johnson.
\newblock {An fMRI investigation of the cultural specificity of music memory}.
\newblock \emph{{Social Cognitive and Affective Neuroscience}}, 5\penalty0
  (2-3):\penalty0 282--291, 2009.
\newblock ISSN 1749-5016.
\newblock \doi{10.1093/scan/nsp048}.

\bibitem[der Nederlanden et~al.(2020)der Nederlanden, Joanisse, and
  Grahn]{Nederlanden2020}
C.~M. V.~B. der Nederlanden, M.~F. Joanisse, and J.~A. Grahn.
\newblock Music as a scaffold for listening to speech: Better neural
  phase-locking to song than speech.
\newblock \emph{{NeuroImage}}, 214:\penalty0 116767, jul 2020.
\newblock \doi{10.1016/j.neuroimage.2020.116767}.

\bibitem[Dumas(2013)]{Dumas2013MelodiesInSpace}
R.~Dumas.
\newblock Melodies in space: neural processing of musical features, Mar. 2013.
\newblock Type: Thesis.

\bibitem[Euler(1774)]{Euler:1774}
L.~Euler.
\newblock De harmoniae veris principiis per speculum musicum repraesentatis.
\newblock \emph{Novi Commentarii academiae scientiarum Petropolitanae},
  18:\penalty0 330--353, 1774.
\newblock URL \url{http://www.math.dartmouth.edu/~euler/pages/E457.html}.

\bibitem[{Euler, Leonhard, 1707-1783}(2017)]{Euler2017}
{Euler, Leonhard, 1707-1783}.
\newblock Tentamen novae theoriae musicae, 2017.

\bibitem[Fechner(1860)]{fechner1860elemente}
G.~Fechner.
\newblock \emph{Elemente der Psychosophysik}.
\newblock Number Bd. 2 in Elemente der Psychosophysik. Breitkopf und Härtel,
  1860.
\newblock URL \url{https://books.google.to/books?id=g3kPAAAAQAAJ}.

\bibitem[Fell and Axmacher(2011)]{Fell2011}
J.~Fell and N.~Axmacher.
\newblock The role of phase synchronization in memory processes.
\newblock \emph{Nature Reviews Neuroscience}, 12\penalty0 (2):\penalty0
  105--118, jan 2011.
\newblock \doi{10.1038/nrn2979}.

\bibitem[Feng(2012)]{feng2012music}
J.~Q. Feng.
\newblock Music in terms of science, 2012.

\bibitem[Gazor and Shoghi(2021)]{Gazor2021}
M.~Gazor and A.~Shoghi.
\newblock Bifurcation control and sound intensities in musical art.
\newblock \emph{Journal of Differential Equations}, 293:\penalty0 86--110, aug
  2021.
\newblock \doi{10.1016/j.jde.2021.05.022}.

\bibitem[Gazor and Shoghi(2022)]{Gazor2022}
M.~Gazor and A.~Shoghi.
\newblock Tone colour in music and bifurcation control.
\newblock \emph{Journal of Differential Equations}, 326:\penalty0 129--163, jul
  2022.
\newblock \doi{10.1016/j.jde.2022.04.011}.

\bibitem[Geer et~al.(1962)Geer, Levelt, and
  Plomp]{Geer1962connotationmusicalconsonance}
J.~V.~D. Geer, W.~Levelt, and R.~Plomp.
\newblock The connotation of musical consonance.
\newblock \emph{Acta Psychologica}, 20:\penalty0 308--319, 1962.
\newblock \doi{10.1016/0001-6918(62)90028-8}.

\bibitem[Genuys(2019)]{Genuys2019}
G.~Genuys.
\newblock Pseudo-distances between chords of different cardinality on
  generalized voice-leading spaces.
\newblock \emph{Journal of Mathematics and Music}, 13\penalty0 (3):\penalty0
  193--206, jul 2019.
\newblock \doi{10.1080/17459737.2019.1622809}.

\bibitem[Gilchrist et~al.(2005)Gilchrist, Jerwood, and Ismaiel]{Gilchrist2005}
J.~M. Gilchrist, D.~Jerwood, and H.~S. Ismaiel.
\newblock Comparing and unifying slope estimates across psychometric function
  models.
\newblock \emph{Perception {\&} Psychophysics}, 67\penalty0 (7):\penalty0
  1289--1303, oct 2005.
\newblock \doi{10.3758/bf03193560}.

\bibitem[Gill and Purves(2009)]{Gill2009}
K.~Z. Gill and D.~Purves.
\newblock A biological rationale for musical scales.
\newblock \emph{{PLoS} {ONE}}, 4\penalty0 (12):\penalty0 e8144, dec 2009.
\newblock \doi{10.1371/journal.pone.0008144}.

\bibitem[Hall and Hess(1984)]{Hall1984}
D.~E. Hall and J.~T. Hess.
\newblock Perception of musical interval tuning.
\newblock \emph{Music Perception}, 2\penalty0 (2):\penalty0 166--195, 1984.
\newblock \doi{10.2307/40285290}.

\bibitem[Harrison(2021)]{Harrison2021}
P.~M.~C. Harrison.
\newblock Three questions concerning consonance~perception.
\newblock \emph{Music Perception}, 38\penalty0 (3):\penalty0 337--339, feb
  2021.
\newblock \doi{10.1525/mp.2021.38.3.337}.

\bibitem[Harrison and
  Pearce(2020{\natexlab{a}})]{Harrison2020ComputationalCognitiveModel}
P.~M.~C. Harrison and M.~T. Pearce.
\newblock {A Computational Cognitive Model for the Analysis and Generation of
  Voice Leadings}.
\newblock \emph{{Music Perception: An Interdisciplinary Journal}}, 37\penalty0
  (3):\penalty0 208--224, 2020{\natexlab{a}}.
\newblock ISSN 0730-7829.
\newblock \doi{10.1525/mp.2020.37.3.208}.
\newblock URL \url{https://mp.ucpress.edu/content/37/3/208}.

\bibitem[Harrison and
  Pearce(2020{\natexlab{b}})]{Harrison2020Simultaneousconsonancemusic}
P.~M.~C. Harrison and M.~T. Pearce.
\newblock Simultaneous consonance in music perception and composition.
\newblock \emph{Psychological Review}, 127\penalty0 (2):\penalty0 216--244, mar
  2020{\natexlab{b}}.
\newblock \doi{10.1037/rev0000169}.

\bibitem[Harrison et~al.(2021)Harrison, Marjieh, and
  Jacoby]{Harrison2021CharacterizingSubjectivePleasantness}
P.~M.~C. Harrison, R.~Marjieh, and N.~Jacoby.
\newblock Characterizing the subjective pleasantness of tone combinations as a
  function of intervallic and spectral structure.
\newblock In \emph{Proceedings of the 16th International Conference on Music
  Perception and Cognition jointly organised with the 11th triennial conference
  of ESCOM}, 2021.

\bibitem[Hartmann(2004)]{Hartmann2004}
W.~M. Hartmann.
\newblock \emph{Signals, Sound, and Sensation}.
\newblock American Inst. of Physics, Sept. 2004.
\newblock ISBN 1563962837.
\newblock URL
  \url{https://www.ebook.de/de/product/1330322/william_m_hartmann_signals_sound_and_sensation.html}.

\bibitem[Hartmann et~al.(1985)Hartmann, Rakerd, and Packard]{Hartmann1985}
W.~M. Hartmann, B.~Rakerd, and T.~N. Packard.
\newblock On measuring the frequency-difference limen for short tones.
\newblock \emph{Perception {\&} Psychophysics}, 38\penalty0 (3):\penalty0
  199--207, may 1985.
\newblock \doi{10.3758/bf03207145}.

\bibitem[Helmholtz(1954)]{Helmholtz1954SensationsTone}
H.~Helmholtz.
\newblock \emph{On the Sensations of Tone}.
\newblock DOVER PUBN INC, June 1954.
\newblock ISBN 0486607534.
\newblock URL
  \url{https://www.ebook.de/de/product/3303046/hermann_helmholtz_on_the_sensations_of_tone.html}.

\bibitem[Hinrichsen(2012)]{Hinrichsen2012}
H.~Hinrichsen.
\newblock Entropy-based tuning of musical instruments.
\newblock \emph{Revista Brasileira de Ensino de F{\'{\i}}sica}, 34\penalty0
  (2):\penalty0 1--8, jun 2012.
\newblock \doi{10.1590/s1806-11172012000200004}.

\bibitem[Hove et~al.(2014)Hove, Marie, Bruce, and Trainor]{Hove2014}
M.~J. Hove, C.~Marie, I.~C. Bruce, and L.~J. Trainor.
\newblock Superior time perception for lower musical pitch explains why
  bass-ranged instruments lay down musical rhythms.
\newblock \emph{Proceedings of the National Academy of Sciences}, 111\penalty0
  (28):\penalty0 10383--10388, jun 2014.
\newblock \doi{10.1073/pnas.1402039111}.

\bibitem[Huckemann et~al.(2010)Huckemann, Hotzand, and Munk]{Huckemann2010}
S.~Huckemann, T.~Hotzand, and A.~Munk.
\newblock Intrinsic shape analysis: Geodesic pca for riemannian manifolds
  modulo isometric lie group actions.
\newblock \emph{Statistica Sinica}, 20, 01 2010.

\bibitem[Hughes(2022)]{Hughes2022}
J.~R. Hughes.
\newblock Generalizing the orbifold model for voice leading.
\newblock \emph{Mathematics}, 10\penalty0 (6):\penalty0 939, mar 2022.
\newblock \doi{10.3390/math10060939}.

\bibitem[Hugo~Fastl(2006)]{HugoFastl2006}
E.~Z. Hugo~Fastl.
\newblock \emph{Psychoacoustics}.
\newblock Springer-Verlag GmbH, Dec. 2006.
\newblock ISBN 3540231595.
\newblock URL
  \url{https://www.ebook.de/de/product/4013268/hugo_fastl_eberhard_zwicker_psychoacoustics.html}.

\bibitem[Huron(2006)]{Huron2006}
D.~Huron.
\newblock \emph{Sweet Anticipation}.
\newblock The {MIT} Press, 2006.
\newblock \doi{10.7551/mitpress/6575.001.0001}.

\bibitem[Johnson(2000)]{Johnson2000}
K.~O. Johnson.
\newblock Neural coding.
\newblock \emph{Neuron}, 26\penalty0 (3):\penalty0 563--566, jun 2000.
\newblock \doi{10.1016/s0896-6273(00)81193-9}.

\bibitem[Jordan(2013)]{jordan2013advancing}
B.~Jordan.
\newblock \emph{Advancing Ethnography in Corporate Environments: Challenges and
  Emerging Opportunities}.
\newblock Left Coast Press, 2013.
\newblock ISBN 9781611322200.
\newblock URL \url{https://books.google.de/books?id=1sk6OEgTVhYC}.

\bibitem[Kac(1966)]{Kac1966}
M.~Kac.
\newblock Can one hear the shape of a drum?
\newblock \emph{The American Mathematical Monthly}, 73\penalty0 (4):\penalty0
  1, apr 1966.
\newblock \doi{10.2307/2313748}.

\bibitem[Kopiez(2003)]{Kopiez2003}
R.~Kopiez.
\newblock Intonation of harmonic intervals: Adaptability of expert musicians to
  equal temperament and just intonation.
\newblock \emph{Music Perception}, 20\penalty0 (4):\penalty0 383--410, 2003.
\newblock \doi{10.1525/mp.2003.20.4.383}.

\bibitem[Krumhansl(1995)]{Krumhansl1995}
C.~L. Krumhansl.
\newblock Music psychology and music theory: Problems and prospects.
\newblock \emph{Music Theory Spectrum}, 17\penalty0 (1):\penalty0 53--80, apr
  1995.
\newblock \doi{10.2307/745764}.

\bibitem[Krumhansl(2001)]{Krumhansl.2001}
C.~L. Krumhansl.
\newblock \emph{{Cognitive Foundations of Musical Pitch}}.
\newblock {Oxford Psychology Series}. {Oxford University Press}, 2001.
\newblock ISBN 9780190287443.
\newblock URL \url{https://books.google.de/books?id=J4dJCAAAQBAJ}.

\bibitem[Lahdelma and Eerola(2015)]{Lahdelma_2015}
I.~Lahdelma and T.~Eerola.
\newblock Theoretical proposals on how vertical harmony may convey nostalgia
  and longing in music.
\newblock \emph{Empirical Musicology Review}, 10\penalty0 (3):\penalty0 245,
  sep 2015.
\newblock \doi{10.18061/emr.v10i3.4534}.

\bibitem[Lahdelma and Eerola(2016)]{Lahdelma_2016}
I.~Lahdelma and T.~Eerola.
\newblock Mild dissonance preferred over consonance in single chord perception.
\newblock \emph{i-Perception}, 7\penalty0 (3):\penalty0 204166951665581, jun
  2016.
\newblock \doi{10.1177/2041669516655812}.

\bibitem[Lahdelma and Eerola(2020)]{Lahdelma2020Culturalfamiliaritymusical}
I.~Lahdelma and T.~Eerola.
\newblock Cultural familiarity and musical expertise impact the pleasantness of
  consonance/dissonance but not its perceived tension.
\newblock \emph{Scientific Reports}, 10\penalty0 (1), may 2020.
\newblock \doi{10.1038/s41598-020-65615-8}.

\bibitem[Lange(2020)]{Lange2020}
C.~Lange.
\newblock Orbifolds from a metric viewpoint.
\newblock \emph{Geometriae Dedicata}, 209\penalty0 (1):\penalty0 43--57, mar
  2020.
\newblock \doi{10.1007/s10711-020-00521-x}.

\bibitem[Langner(1997)]{Langner1997Temporalprocessingpitch}
G.~Langner.
\newblock Temporal processing of pitch in the auditory system.
\newblock \emph{Journal of New Music Research}, 26\penalty0 (2):\penalty0
  116--132, jun 1997.
\newblock \doi{10.1080/09298219708570721}.

\bibitem[Langner and Benson(2015)]{Langner2015NeuralCodePitch}
G.~Langner and C.~Benson.
\newblock \emph{The Neural Code of Pitch and Harmony}.
\newblock Cambridge University Press, 2015.
\newblock \doi{10.1017/cbo9781139050852}.

\bibitem[Large(2010)]{Large2010}
E.~W. Large.
\newblock A dynamical systems approach to musical tonality.
\newblock In \emph{Nonlinear Dynamics in Human Behavior}, pages 193--211.
  Springer Berlin Heidelberg, 2010.
\newblock \doi{10.1007/978-3-642-16262-6_9}.

\bibitem[Leino et~al.(2007)Leino, Brattico, Tervaniemi, and Vuust]{Leino2007}
S.~Leino, E.~Brattico, M.~Tervaniemi, and P.~Vuust.
\newblock Representation of harmony rules in the human brain: Further evidence
  from event-related potentials.
\newblock \emph{Brain Research}, 1142:\penalty0 169--177, apr 2007.
\newblock \doi{10.1016/j.brainres.2007.01.049}.

\bibitem[Leman(2000)]{Leman2000}
M.~Leman.
\newblock Visualization and calculation of the roughness ofacoustical music
  signals using the synchronization index model.
\newblock In \emph{Proceedingsof the COSTG-6 Conference on Digital Audio
  Effects (DAFX-00)}, Verona, Italy, 2000.

\bibitem[Lerdahl(2001)]{Lerdahl.2001}
F.~Lerdahl.
\newblock \emph{{Tonal Pitch Space}}.
\newblock {Oxford University Press}, New York, 2001.

\bibitem[Lodish(2000)]{Lodish2000}
H.~Lodish.
\newblock \emph{Molecular cell biology}.
\newblock W.H. Freeman, New York, 2000.
\newblock ISBN 0716731363.

\bibitem[Lowet et~al.(2016)Lowet, Roberts, Bonizzi, Karel, and
  Weerd]{Lowet2016}
E.~Lowet, M.~J. Roberts, P.~Bonizzi, J.~Karel, and P.~D. Weerd.
\newblock Quantifying neural oscillatory synchronization: A comparison between
  spectral coherence and phase-locking value approaches.
\newblock \emph{{PLOS} {ONE}}, 11\penalty0 (1):\penalty0 e0146443, jan 2016.
\newblock \doi{10.1371/journal.pone.0146443}.

\bibitem[Macmillan and Creelman(2004)]{Macmillan2004}
N.~A. Macmillan and C.~D. Creelman.
\newblock \emph{Detection Theory}.
\newblock Psychology Press, sep 2004.
\newblock \doi{10.4324/9781410611147}.

\bibitem[Maher(1976)]{Maher1976NeedResolutionRatings}
T.~F. Maher.
\newblock "need for resolution" ratings for harmonic musical intervals.
\newblock \emph{Journal of Cross-Cultural Psychology}, 7\penalty0 (3):\penalty0
  259--276, sep 1976.
\newblock \doi{10.1177/002202217673001}.

\bibitem[Marcus(1993)]{Marcus1993}
S.~Marcus.
\newblock The interface between theory and practice: Intonation in arab music.
\newblock \emph{Asian Music}, 24\penalty0 (2):\penalty0 39, 1993.
\newblock \doi{10.2307/834466}.

\bibitem[Marjieh et~al.(2022)Marjieh, Harrison, Lee, Deligiannaki, and
  Jacoby]{Marjieh2022}
R.~Marjieh, P.~M.~C. Harrison, H.~Lee, F.~Deligiannaki, and N.~Jacoby.
\newblock Reshaping musical consonance with timbral manipulations and massive
  online experiments.
\newblock \emph{bioRxiv}, 2022.
\newblock \doi{10.1101/2022.06.14.496070}.
\newblock URL
  \url{https://www.biorxiv.org/content/early/2022/06/17/2022.06.14.496070}.

\bibitem[Mattson(2014)]{Mattson.2014}
M.~P. Mattson.
\newblock {Superior pattern processing is the essence of the evolved human
  brain}.
\newblock \emph{{Frontiers in neuroscience}}, 8:\penalty0 265, 2014.
\newblock ISSN 1662-453X.
\newblock \doi{10.3389/fnins.2014.00265}.

\bibitem[Michor(2008)]{Michor2008}
P.~W. Michor.
\newblock \emph{Topics in differential geometry}.
\newblock American Mathematical Society, 2008.
\newblock ISBN 9780821820032.

\bibitem[Milne and Holland(2016)]{Milne2016}
A.~J. Milne and S.~Holland.
\newblock Empirically testing tonnetz, voice-leading, and spectral models of
  perceived triadic distance.
\newblock \emph{Journal of Mathematics and Music}, 10\penalty0 (1):\penalty0
  59--85, jan 2016.
\newblock \doi{10.1080/17459737.2016.1152517}.

\bibitem[{Milne, Andrew J.} et~al.(2011){Milne, Andrew J.}, {Sethares, William
  A.}, Laney, and {Sharp, David B.}]{MilneAndrewJ..2011}
{Milne, Andrew J.}, {Sethares, William A.}, R.~Laney, and {Sharp, David B.}
\newblock {Modelling the similarity of pitch collections with expectation
  tensors}.
\newblock \emph{{Journal of Mathematics and Music}}, 5\penalty0 (1):\penalty0
  1--20, 2011.
\newblock \doi{10.1080/17459737.2011.573678}.

\bibitem[Moore et~al.(1984)Moore, Glasberg, and
  Shailer]{Moore1984Frequencyintensitydifference}
B.~C.~J. Moore, B.~R. Glasberg, and M.~J. Shailer.
\newblock Frequency and intensity difference limens for harmonics within
  complex tones.
\newblock \emph{The Journal of the Acoustical Society of America}, 75\penalty0
  (2):\penalty0 550--561, feb 1984.
\newblock \doi{10.1121/1.390527}.

\bibitem[Moore et~al.(1985)Moore, Peters, and Glasberg]{Moore1985}
B.~C.~J. Moore, R.~W. Peters, and B.~R. Glasberg.
\newblock Thresholds for the detection of inharmonicity in complex tones.
\newblock \emph{The Journal of the Acoustical Society of America}, 77\penalty0
  (5):\penalty0 1861--1867, may 1985.
\newblock \doi{10.1121/1.391937}.

\bibitem[Moore et~al.(1986)Moore, Glasberg, and Peters]{Moore1986}
B.~C.~J. Moore, B.~R. Glasberg, and R.~W. Peters.
\newblock Thresholds for hearing mistuned partials as separate tones in
  harmonic complexes.
\newblock \emph{The Journal of the Acoustical Society of America}, 80\penalty0
  (2):\penalty0 479--483, aug 1986.
\newblock \doi{10.1121/1.394043}.

\bibitem[Owen(2000)]{Owen.2000}
H.~Owen.
\newblock \emph{{Music Theory Resource Book}}.
\newblock {Oxford and New York: Oxford University Press}, 2000.
\newblock ISBN 0-19-511539-2.

\bibitem[Pag{\`{e}}s-Portabella and Toro(2019)]{PagesPortabella2019}
C.~Pag{\`{e}}s-Portabella and J.~M. Toro.
\newblock Dissonant endings of chord progressions elicit a larger {ERAN} than
  ambiguous endings in musicians.
\newblock \emph{Psychophysiology}, 57\penalty0 (2), sep 2019.
\newblock \doi{10.1111/psyp.13476}.

\bibitem[Parncutt and Hair(2012)]{Parncutt2012Consonancedissonancemusic}
R.~Parncutt and G.~Hair.
\newblock Consonance and dissonance in music theory and psychology:
  Disentangling dissonant dichotomies.
\newblock \emph{JOurnal of interdisciplinary music studies}, 2012.
\newblock \doi{10.4407/jims.2011.11.002}.

\bibitem[Partridge and Partridge(2003)]{Partridge2003}
L.~D. Partridge and L.~D. Partridge.
\newblock From reception to pattern recognition and perception.
\newblock In \emph{Nervous System Actions and Interactions}, pages 145--174.
  Springer {US}, 2003.
\newblock \doi{10.1007/978-1-4615-0425-2_8}.

\bibitem[Pearce and Wiggins(2012)]{Pearce2012}
M.~T. Pearce and G.~A. Wiggins.
\newblock Auditory expectation: The information dynamics of music perception
  and cognition.
\newblock \emph{Topics in Cognitive Science}, 4\penalty0 (4):\penalty0
  625--652, jul 2012.
\newblock \doi{10.1111/j.1756-8765.2012.01214.x}.

\bibitem[Pflaum(2003)]{Pflaum2003AnalyticGeometricStudy}
M.~Pflaum.
\newblock \emph{Analytic and Geometric Study of Stratified Spaces:
  Contributions to Analytic and Geometric Aspects}.
\newblock Lecture Notes in Mathematics. Springer Berlin Heidelberg, 2003.
\newblock ISBN 9783540454366.

\bibitem[Railsback(1938)]{Railsback1938ScaleTemperamentas}
O.~L. Railsback.
\newblock Scale temperament as applied to piano tuning.
\newblock \emph{The Journal of the Acoustical Society of America}, 9\penalty0
  (3):\penalty0 274--274, jan 1938.
\newblock \doi{10.1121/1.1902056}.

\bibitem[Randall and Khan(2010)]{Randall2010}
R.~R. Randall and B.~Khan.
\newblock {Lerdahl's tonal pitch space model and associated metric spaces}.
\newblock \emph{{Journal of Mathematics and Music}}, 4\penalty0 (3):\penalty0
  121--131, 2010.
\newblock \doi{10.1080/17459737.2010.529654}.

\bibitem[Randel(1999)]{Randel.1999}
D.~M. Randel.
\newblock \emph{{The Harvard Concise Dictionary of Music and Musicians}}.
\newblock {Belknap Press}, 1999.
\newblock ISBN 0-674-00084-6.

\bibitem[Ratcliffe(2007)]{Ratcliffe2007}
J.~Ratcliffe.
\newblock \emph{Foundations of Hyperbolic Manifolds}.
\newblock Springer New York, 2007.
\newblock \doi{10.1007/978-0-387-47322-2}.

\bibitem[Rickles(2016)]{Rickles.2016}
D.~Rickles.
\newblock {Spaces}.
\newblock In P.~Humphreys, editor, \emph{{The Oxford Handbook of Philosophy of
  Science}}, Oxford, United Kingdom, 2016. {Oxford University Press}.
\newblock URL
  \url{https://www.oxfordhandbooks.com/view/10.1093/oxfordhb/9780199368815.001.0001/oxfordhb-9780199368815-e-31}.

\bibitem[Roederer(2008)]{Roederer2008}
J.~G. Roederer.
\newblock \emph{The Physics and Psychophysics of Music}.
\newblock Springer-Verlag GmbH, Dec. 2008.
\newblock ISBN 9780387094748.
\newblock URL
  \url{https://www.ebook.de/de/product/12470007/juan_g_roederer_the_physics_and_psychophysics_of_music.html}.

\bibitem[Rohrmeier(2013)]{Rohrmeier2013}
M.~Rohrmeier.
\newblock Musical expectancy. bridging music theory, cognitive and
  computational approaches.
\newblock \emph{Zeitschrift der Gesellschaft für Musiktheorie [Journal of the
  German-Speaking Society of Music Theory]}, 10\penalty0 (2):\penalty0
  343--371, 2013.
\newblock \doi{10.31751/724}.

\bibitem[Rossing(2002)]{Rossing2002}
T.~Rossing.
\newblock \emph{The {S}cience of {S}ound}.
\newblock Addison Wesley, San Francisco, 2002.
\newblock ISBN 9780805385656.

\bibitem[Sauv{\'{e}} et~al.(2021)Sauv{\'{e}}, Cho, and Zendel]{Sauve2021}
S.~A. Sauv{\'{e}}, A.~Cho, and B.~R. Zendel.
\newblock Mapping tonal hierarchy in the brain.
\newblock \emph{Neuroscience}, 465:\penalty0 187--202, jun 2021.
\newblock \doi{10.1016/j.neuroscience.2021.03.019}.

\bibitem[Schmuckler(1989)]{Schmuckler1989}
M.~A. Schmuckler.
\newblock Expectation in music: Investigation of melodic and harmonic
  processes.
\newblock \emph{Music Perception}, 7\penalty0 (2):\penalty0 109--149, 1989.
\newblock \doi{10.2307/40285454}.

\bibitem[Sch{\"o}nberg et~al.(2000)Sch{\"o}nberg, Black, and
  Stein]{Schonberg.2000}
A.~Sch{\"o}nberg, L.~Black, and L.~Stein.
\newblock \emph{{Style and idea: Selected writings of Arnold Sch{\"o}nberg}}.
\newblock {University of California Press}, Berkley, 2000.
\newblock ISBN 0-520-05294-3.

\bibitem[Schwitzgebel and White(2021)]{Schwitzgebel2021}
E.~Schwitzgebel and C.~W. White.
\newblock Effects of chord inversion and bass patterns on harmonic expectancy
  in musicians.
\newblock \emph{Music Perception}, 39\penalty0 (1):\penalty0 41--62, sep 2021.
\newblock \doi{10.1525/mp.2021.39.1.41}.

\bibitem[Seger et~al.(2013)Seger, Spiering, Sares, Quraini, Alpeter, David, and
  Thaut]{Seger2013}
C.~A. Seger, B.~J. Spiering, A.~G. Sares, S.~I. Quraini, C.~Alpeter, J.~David,
  and M.~H. Thaut.
\newblock Corticostriatal contributions to musical expectancy perception.
\newblock \emph{Journal of Cognitive Neuroscience}, 25\penalty0 (7):\penalty0
  1062--1077, jul 2013.
\newblock \doi{10.1162/jocn_a_00371}.

\bibitem[Sethares(2005)]{sethares2005tuning}
W.~Sethares.
\newblock \emph{Tuning, Timbre, Spectrum, Scale}.
\newblock Springer London, 2005.
\newblock ISBN 9781852337971.
\newblock URL \url{https://books.google.de/books?id=KChoKKhjOb0C}.

\bibitem[Sethares(1993)]{Sethares1993}
W.~A. Sethares.
\newblock Local consonance and the relationship between timbre and scale.
\newblock \emph{The Journal of the Acoustical Society of America}, 94\penalty0
  (3):\penalty0 1218--1228, sep 1993.
\newblock \doi{10.1121/1.408175}.

\bibitem[Shepard(1964)]{Shepard.1964}
R.~N. Shepard.
\newblock {Circularity in judgments of relative pitch}.
\newblock \emph{{J. Acoust. Soc. Am.}}, 36:\penalty0 2345--2353, 1964.

\bibitem[Sinz et~al.(2020)Sinz, Sachgau, Henninger, Benda, and Grewe]{Sinz2020}
F.~H. Sinz, C.~Sachgau, J.~Henninger, J.~Benda, and J.~Grewe.
\newblock Simultaneous spike-time locking to multiple frequencies.
\newblock \emph{Journal of Neurophysiology}, 123\penalty0 (6):\penalty0
  2355--2372, jun 2020.
\newblock \doi{10.1152/jn.00615.2019}.

\bibitem[Smith and Abel(1999)]{Smith1999}
J.~Smith and J.~Abel.
\newblock Bark and {ERB} bilinear transforms.
\newblock \emph{{IEEE} Transactions on Speech and Audio Processing}, 7\penalty0
  (6):\penalty0 697--708, 1999.
\newblock \doi{10.1109/89.799695}.

\bibitem[Stefano and
  Bertolaso(2014)]{Stefano2014UnderstandingMusicalConsonance}
N.~D. Stefano and M.~Bertolaso.
\newblock Understanding musical consonance and dissonance: Epistemological
  considerations from a systemic perspective.
\newblock \emph{Systems}, 2\penalty0 (4):\penalty0 566--575, oct 2014.
\newblock \doi{10.3390/systems2040566}.

\bibitem[Stolzenburg(2015)]{Stolzenburg2015Harmonyperceptionperiodicity}
F.~Stolzenburg.
\newblock Harmony perception by periodicity detection.
\newblock \emph{Journal of Mathematics and Music}, 9\penalty0 (3):\penalty0
  215--238, aug 2015.
\newblock \doi{10.1080/17459737.2015.1033024}.

\bibitem[Stumpf(2013{\natexlab{a}})]{Stumpf2013_1}
C.~Stumpf.
\newblock \emph{Tonpsychologie 1}.
\newblock Cambridge University Press, May 2013{\natexlab{a}}.
\newblock ISBN 110806177X.
\newblock URL
  \url{https://www.ebook.de/de/product/20825970/carl_stumpf_tonpsychologie.html}.

\bibitem[Stumpf(2013{\natexlab{b}})]{Stumpf2013_2}
C.~Stumpf.
\newblock \emph{Tonpsychologie 2}.
\newblock Cambridge University Press, May 2013{\natexlab{b}}.
\newblock ISBN 1108061788.
\newblock URL
  \url{https://www.ebook.de/de/product/20825969/carl_stumpf_tonpsychologie.html}.

\bibitem[Sutcliffe(2011)]{Sutcliffe.2011}
T.~Sutcliffe.
\newblock \emph{{Syntactic Structures in Music}}.
\newblock http://www.harmony.org.uk/, 2011.
\newblock URL \url{https://www.harmony.org.uk}.

\bibitem[Thanwerdas(2022)]{Thanwerdas2022RiemannianStratifiedGeometries}
Y.~Thanwerdas.
\newblock \emph{{Riemannian and stratified geometries on covariance and
  correlation matrices}}.
\newblock Theses, {Universit{\'e} C{\^o}te d'Azur}, May 2022.
\newblock URL \url{https://hal.archives-ouvertes.fr/tel-03698752}.

\bibitem[Thurston(1997)]{Thurston1997}
W.~P. Thurston.
\newblock \emph{Three-dimensional geometry and topology}.
\newblock Princeton University Press, 1997.
\newblock ISBN 0691083045.

\bibitem[TRAMO(2005)]{TRAMO_2005}
M.~J. TRAMO.
\newblock Neurophysiology and neuroanatomy of pitch perception: Auditory
  cortex.
\newblock \emph{Annals of the New York Academy of Sciences}, 1060\penalty0
  (1):\penalty0 148--174, dec 2005.
\newblock \doi{10.1196/annals.1360.011}.

\bibitem[Tramo et~al.(2001)Tramo, Cariani, Delgutte, and Braida]{Tramo.2001}
M.~J. Tramo, P.~A. Cariani, B.~Delgutte, and L.~D. Braida.
\newblock {Neurobiological foundations for the theory of harmony in western
  tonal music}.
\newblock \emph{{Annals of the New York Academy of Sciences}}, 930:\penalty0
  92--116, 2001.
\newblock ISSN 0077-8923.
\newblock \doi{10.1111/j.1749-6632.2001.tb05727.x}.

\bibitem[Tymoczko(2006)]{Tymoczko.2006}
D.~Tymoczko.
\newblock {The Geometry of Musical Chords}.
\newblock \emph{{Science}}, 313:\penalty0 72--74, 2006.

\bibitem[Tymoczko(2009)]{Tymoczko.2009}
D.~Tymoczko.
\newblock {Three Conceptions of Musical Distance}.
\newblock In \emph{{Mathematics and Computation in Music, eds. Elaine Chew,
  Adrian Childs, and Ching-Hua Chuan}}, pages 258--273. Springer, Heidelberg,
  2009.

\bibitem[Tymoczko(2011)]{Tymoczko.2011}
D.~Tymoczko.
\newblock \emph{{A Geometry of Music: Harmony and Counterpoint in the Extended
  Common Practice}}.
\newblock {Oxford Studies in Music Theory}. {Oxford University Press}, 2011.
\newblock ISBN 9780199887507.
\newblock URL \url{https://books.google.de/books?id=ODSt58Yk2YYC}.

\bibitem[Tymoczko(2012)]{Tymoczko_2012}
D.~Tymoczko.
\newblock The generalized tonnetz.
\newblock \emph{Journal of Music Theory}, 56\penalty0 (1):\penalty0 1--52, mar
  2012.
\newblock \doi{10.1215/00222909-1546958}.

\bibitem[Valla et~al.(2017)Valla, Alappatt, Mathur, and Singh]{Valla2017}
J.~M. Valla, J.~A. Alappatt, A.~Mathur, and N.~C. Singh.
\newblock Music and emotion{\textemdash}a case for north indian classical
  music.
\newblock \emph{Frontiers in Psychology}, 8, dec 2017.
\newblock \doi{10.3389/fpsyg.2017.02115}.

\bibitem[Varshney and Sun(2013)]{Varshney2013LogarithmicPerception}
L.~R. Varshney and J.~Z. Sun.
\newblock Why do we perceive logarithmically?
\newblock \emph{Significance}, 10\penalty0 (1):\penalty0 28--31, 2013.
\newblock \doi{https://doi.org/10.1111/j.1740-9713.2013.00636.x}.
\newblock URL
  \url{https://rss.onlinelibrary.wiley.com/doi/abs/10.1111/j.1740-9713.2013.00636.x}.

\bibitem[Vencovsk{\'{y}}(2016)]{Vencovsky2016}
V.~Vencovsk{\'{y}}.
\newblock Roughness prediction based on a model of cochlear hydrodynamics.
\newblock \emph{Archives of Acoustics}, 41\penalty0 (2):\penalty0 189--201, jun
  2016.
\newblock \doi{10.1515/aoa-2016-0019}.

\bibitem[Vos(1986)]{Vos1986}
J.~Vos.
\newblock Purity ratings of tempered fifths and major thirds.
\newblock \emph{Music Perception}, 3\penalty0 (3):\penalty0 221--257, 1986.
\newblock \doi{10.2307/40285335}.

\bibitem[Vuust et~al.(2022)Vuust, Heggli, Friston, and Kringelbach]{Vuust2022}
P.~Vuust, O.~A. Heggli, K.~J. Friston, and M.~L. Kringelbach.
\newblock Music in the brain.
\newblock \emph{Nature Reviews Neuroscience}, 23\penalty0 (5):\penalty0
  287--305, mar 2022.
\newblock \doi{10.1038/s41583-022-00578-5}.

\bibitem[Wall et~al.(2020)Wall, Lieck, Neuwirth, and Rohrmeier]{Wall2020}
L.~Wall, R.~Lieck, M.~Neuwirth, and M.~Rohrmeier.
\newblock The impact of voice leading and harmony on musical expectancy.
\newblock \emph{Scientific Reports}, 10\penalty0 (1), apr 2020.
\newblock \doi{10.1038/s41598-020-61645-4}.

\bibitem[Wichmann and Hill(2001{\natexlab{a}})]{Wichmann2001}
F.~A. Wichmann and N.~J. Hill.
\newblock The psychometric function: I. fitting, sampling, and goodness of fit.
\newblock \emph{Perception {\&} Psychophysics}, 63\penalty0 (8):\penalty0
  1293--1313, nov 2001{\natexlab{a}}.
\newblock \doi{10.3758/bf03194544}.

\bibitem[Wichmann and Hill(2001{\natexlab{b}})]{Wichmann2001a}
F.~A. Wichmann and N.~J. Hill.
\newblock The psychometric function: {II}. bootstrap-based confidence intervals
  and sampling.
\newblock \emph{Perception {\&} Psychophysics}, 63\penalty0 (8):\penalty0
  1314--1329, nov 2001{\natexlab{b}}.
\newblock \doi{10.3758/bf03194545}.

\bibitem[Wilkerson(2014)]{wilkerson2014harmony}
D.~S. Wilkerson.
\newblock Harmony explained: Progress towards a scientific theory of music,
  2014.

\bibitem[Wright(2009)]{Wright2009}
D.~Wright.
\newblock \emph{Mathematics and music}.
\newblock American Mathematical Society, Providence, R.I, 2009.
\newblock ISBN 9780821848739.

\bibitem[Zhang et~al.(2018)Zhang, Zhou, Chang, and Yang]{Zhang2018}
J.~Zhang, X.~Zhou, R.~Chang, and Y.~Yang.
\newblock Effects of global and local contexts on chord processing: An {ERP}
  study.
\newblock \emph{Neuropsychologia}, 109:\penalty0 149--154, jan 2018.
\newblock \doi{10.1016/j.neuropsychologia.2017.12.016}.

\bibitem[Zwicker et~al.(1957)Zwicker, Flottorp, and Stevens]{Zwicker1957}
E.~Zwicker, G.~Flottorp, and S.~S. Stevens.
\newblock Critical band width in loudness summation.
\newblock \emph{The Journal of the Acoustical Society of America}, 29\penalty0
  (5):\penalty0 548--557, may 1957.
\newblock \doi{10.1121/1.1908963}.

\end{thebibliography}
	
\end{document}